\documentclass[11pt]{article}
\usepackage[centering,includeheadfoot,margin=2 cm]{geometry}
\usepackage{graphics,enumitem,epsfig}
\usepackage{amsfonts,amsmath,amssymb,euscript,color}

\begin{document}

\def\ALERT#1{{\large\bf $\clubsuit$#1$\clubsuit$}}

\numberwithin{equation}{section}
\newtheorem{defin}{Definition}
\newtheorem{theorem}{Theorem}
\newtheorem{notice}{Notice}
\newtheorem{lemma}{Lemma}
\newtheorem{cor}{Corollary}
\newtheorem{example}{Example}
\newtheorem{remark}{Remark}
\newtheorem{conj}{Conjecture}
\def\begproof{\noindent{\bf Proof: }}
\def\endproof{\par\rightline{\vrule height5pt width5pt depth0pt}\medskip}
\def\div{\nabla\cdot}
\def\rot{\nabla\times}
\def\sign{{\rm sign}}
\def\arsinh{{\rm arsinh}}
\def\arcosh{{\rm arcosh}}
\def\diag{{\rm diag}}
\def\const{{\rm const}}
\def\d{\,\mathrm{d}}
\def\eps{\varepsilon}
\def\phi{\varphi}
\def\theta{\vartheta}
\def\N{\mathbb{N}}
\def\R{\mathbb{R}}
\def\C{\hbox{\rlap{\kern.24em\raise.1ex\hbox
      {\vrule height1.3ex width.9pt}}C}}
\def\P{\hbox{\rlap{I}\kern.16em P}}
\def\Q{\hbox{\rlap{\kern.24em\raise.1ex\hbox
      {\vrule height1.3ex width.9pt}}Q}}
\def\M{\hbox{\rlap{I}\kern.16em\rlap{I}M}}
\def\Z{\hbox{\rlap{Z}\kern.20em Z}}
\def\({\begin{eqnarray}}
\def\){\end{eqnarray}}
\def\[{\begin{eqnarray*}}
\def\]{\end{eqnarray*}}
\def\part#1#2{\frac{\partial #1}{\partial #2}}
\def\partk#1#2#3{\frac{\partial^#3 #1}{\partial #2^#3}} 
\def\mat#1{{D #1\over Dt}}
\def\dx{\nabla_x}
\def\dv{\nabla_v}
\def\grad{\nabla}
\def\Norm#1{\left\| #1 \right\|}
\def\pmb#1{\setbox0=\hbox{$#1$}
  \kern-.025em\copy0\kern-\wd0
  \kern-.05em\copy0\kern-\wd0
  \kern-.025em\raise.0433em\box0 }
\def\bar{\overline}
\def\lbar{\underline}
\def\fref#1{(\ref{#1})}
\def\half{\frac{1}{2}}
\def\oo#1{\frac{1}{#1}}

\def\tot#1#2{\frac{\d #1}{\d #2}} 
\def\laplace{\Delta}
\def\d{\,\mathrm{d}}
\def\N{\mathbb{N}}
\def\R{\mathbb{R}}
\def\supp{\mbox{supp }}
\def\eps{\varepsilon}
\def\phi{\varphi}

\def\L{\mathcal{L}}
\def\H{\mathcal{H}}
\def\O{\mathcal{O}}
\def\K{\mathcal{K}}

\def\w{w}
\def\ti#1{{{\tilde{#1}^i}}}
\def\setv#1{\{\pm1\}^{#1}}

\def\comment#1{\textcolor{blue}{\bf [Comment: #1]}}

\newcommand{\boldx}{{\mathbf x}}
\newcommand{\boldv}{{\mathbf v}}
\newcommand{\boldy}{{\mathbf y}}
\newcommand{\boldz}{{\mathbf z}}


\centerline{{\huge From individual to collective behaviour of coupled velocity}}

\vskip 0.5mm

\centerline{{\huge jump processes: a locust example}}


\vskip 5mm

\centerline{
{\large Radek Erban}\footnote{Mathematical Institute, 
University of Oxford, 24-29 St. Giles', Oxford, OX1 3LB, United Kingdom;\\
e-mail: {\it erban@maths.ox.ac.uk}}
\qquad
{\large Jan Ha{\v s}kovec}\footnote{Johann Radon Institute for Computational 
and Applied Mathematics (RICAM), Austrian Academy of Sciences, \\
Altenbergerstra\ss{}e 69, A-4040 Linz, Austria; 
e-mail: {\it jan.haskovec@oeaw.ac.at}}
}

\vskip 6mm


\noindent{\bf Abstract.}
A class of stochastic individual-based models, written in terms of coupled 
velocity jump processes, is presented and analysed. This modelling approach 
incorporates recent experimental findings on behaviour of locusts. 
It exhibits nontrivial dynamics with a ``phase change" behaviour and recovers the 
observed group directional switching. Estimates of the expected switching 
times, in terms of number of individuals and values of the model coefficients,
are obtained using the corresponding Fokker-Planck equation. In the limit
of large populations, a system of two kinetic equations with nonlocal and 
nonlinear right hand side is derived and analyzed. The existence of its 
solutions is proven and the system's long-time behaviour is investigated.
Finally, a first step towards the mean field limit of topological interactions
is made by studying the effect of shrinking the interaction radius in the
individual-based model when the number of individuals grows.

\vskip 5mm


\section{Introduction}\label{sec:Introduction}

Individual-based behaviour in biology can be often modelled 
as a velocity jump process \cite{Othmer:1988:MDB}. Here, the
velocity of an individual is subject to sudden changes (``jumps") 
at random instants. If the velocity changes are 
completely random, this process simply leads to diffusive 
spreading of individuals in an appropriate limit \cite{Hillen:2000:DLT}. 
The situation becomes more complicated whenever the velocity
changes are biased according to an individual's environment.
A classical example is bacterial chemotaxis \cite{Erban:2004:ICB}. 
Individual bacteria change their frequency of velocity changes 
according to their environment. If they swim in a favourable 
direction (e.g. towards a nutrient source), they are less 
likely to change their direction. On the other hand, they are 
more likely to turn if they are heading away from a foodstuff
\cite{Erban:2005:STS}.

In this paper, we modify the velocity jump methodology to model
the behaviour of locusts. Our model is motivated by the recent
experiments of Buhl et al \cite{Buhl:2006:DOM}. They studied an
experimental setting, in which locust nymphs marched in 
a ring-shaped arena. The collective behaviour depended 
strongly on locust density. At low densities, there was a low 
incidence of alignment among individuals. Intermediate densities 
were characterized by long periods of collective motion in
one direction along the arena interrupted by rapid changes of 
group direction. If the density of locusts was further increased, 
the group quickly adopted a common and persistent rotational 
direction. Yates et al \cite{Yates:2009:INC} analysed 
experimental data of Buhl et al \cite{Buhl:2006:DOM} and proposed 
that the frequency of random changes in the direction of
an individual increases when the individual looses the alignment with the
rest of the group. In this paper, we incorporate this observation
into a stochastic individual-based model formulated as a velocity 
jump process. We show that this model, although phenomenologically 
very simple, has the same predictive power as other modelling
approaches previously used in this area \cite{Czirok:1999:CMS,Buhl:2006:DOM}.
In particular, it exhibits (i) a rapid transition from disordered 
movement of individuals to highly aligned collective motion as 
the size of the group grows, and (ii) sudden and rapid switching 
of the group direction, with frequency decreasing as 
group size increases.

The individual-based model is introduced in Section~\ref{sec:Model}. 
The ring-shaped arena, used in Buhl's experiments \cite{Buhl:2006:DOM},
is modelled as one-dimensional interval with periodic boundary
conditions. Locusts march with a constant speed and each individual 
switches its direction randomly. The individual switching frequency 
increases in response to a loss of alignment. In Section~\ref{sec:Fokker-Planck},
the corresponding Fokker-Planck equation is derived for the 
system with global interactions and possible types of qualitative behaviour 
of the system are classified. For the case of ordered group motion, 
where two distinct metastable states exist, an approximate analytic 
formula for the mean switching time between these two states is derived.
Then, in Section~\ref{sec:Kinetic}, the kinetic formulation of the model 
is obtained in the limit as the number of locusts tends to infinity.
The existence of solutions of the kinetic model is shown in Section~\ref{sec:Existence}
and the long time behaviour is investigated in Section~\ref{sec:Long time}.
We conclude with analysis of the dependence of collective behaviour on the
size of the interaction radius of the individuals in Section~\ref{sec:Shrinking}.

\section{Individual based model}\label{sec:Model}
We consider a group of $N$ agents (locusts)
with time-dependent positions $x_i(t)$ and velocities $v_i(t)$, 
$i = 1,\dots, N$. To mimic the ring-shaped arena set-up 
of~\cite{Buhl:2006:DOM}, we assume that the agents move along 
a one-dimensional circle, which we identify with the interval 
$\Omega = [0,1)$ with periodic boundary conditions,
and move either to the right or to the left
with the same unit speed, i.e. 
\(  \label{model0}
x_i(t) \in \Omega, \qquad 
v_i \in \{-1, 1\}
\qquad \mbox{and} \qquad 
\frac{\mbox{d}x_i}{\mbox{d}t} (t) 
= v_i(t).
\)
We define the local average velocity of the ensemble, seen by the 
$i$-th agent, as
\begin{equation}  
   u_i^{loc} = 
   \frac{\sum_{m=1}^N \w(|x_i-x_m|) v_m}{\sum_{m=1}^N \w(|x_i-x_m|)}  \,,
\label{u_i^loc}
\end{equation}
where $\w$ is a weight function defined on $\Omega$ with the properties:

\leftskip 1cm

\begin{enumerate}[leftmargin = 1.6 cm, label = {[}A{}\arabic*{]} ]
\item\label{A1}
$\w$ is bounded and nonnegative on $\Omega$,
\item\label{A2}
$w(0) > 0$.
\end{enumerate}

\par\leftskip 0cm

\noindent
For example, $\w(s) = \chi_{[0,\sigma]}(s)$,
where $\chi_{[0,\sigma]}$ is the characteristic function of the interval $[0,\sigma]$
and $\sigma > 0$ is a interaction radius, satisfies conditions 
\ref{A1} and \ref{A2}. This is a common
choice of $\w$ in biological applications 
\cite{Czirok:1999:CMS, Yates:2009:INC}.
It is worth noting that, due to the assumption~\ref{A2}, 
the definition~\fref{u_i^loc} always makes sense and 
$u_i^{loc} \in [-1,1].$

The agents switch their velocities to the opposite direction
(i.e., from $v_i = 1$ to $v_i = -1$ and vice versa)
based on $N$ independent Poisson processes
with the rates
\[
    \gamma_i = \gamma_0 + b \, \xi(v_i-u_i^{loc}) \,,\qquad i=1,\dots,N\,,
\]
where $\gamma_0 \geq 0$ and $b \geq 0$ are fixed parameters
and the ``response to disalignment'' function $\xi: [-2,2] \to [0,\infty)$ 
is assumed to be convex, differentiable and symmetric with respect to the origin.
Taking the Taylor expansion of $\xi(s)$ around $s=0$, we obtain
\[
   \xi(s) = \alpha_0 + \alpha_2 s^2 + \O(s^3) \,.
\]
We can set $\alpha_0=0$ without loss of generality, because it can be absorbed in
$\gamma_0$.
Since the individuals switch their velocities less frequently when they are aligned
(\cite{Yates:2009:INC}), $\xi(s)$ has a global minimum at $s=0$.
This implies that $\alpha_2 \geq 0$. If $\alpha_2 > 0$, we can set
$\alpha_2 = 1$ by choosing an appropriate time scale.
Therefore, $\xi(s)$ has the general form $s^2 + \O(s^3)$.
For the rest of the paper, we choose the form $\xi(s) = s^2$ for simplicity.
Other choices are certainly possible, for example, in the limiting case $\alpha_2 = 0$
the leading order approximation is given by a higher order term, which, however, complicates
the analysis. However, it is worth noting that the derivation
of the kinetic equation performed in Section~\ref{sec:Kinetic}
is possible, for example, also for $\xi(s) = |s|^n$, $n \geq 3$.

\begin{figure}
\begin{center}
\epsfig{figure=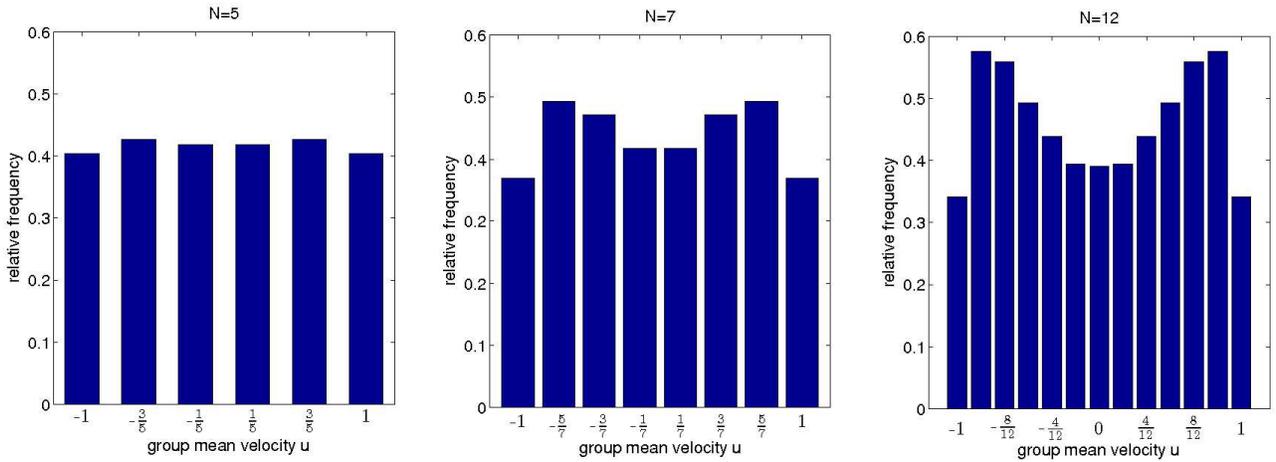, width=1.0\linewidth}
\end{center}
\caption{{\it An example of the transition to ordered motion as $N$ grows.
We use~$\fref{model0}$--$\fref{gammaLR}$ with $b=1$, $\gamma_0=0.2$ and $\w = \chi_{[0,0.2]}$.
Shown are the normalized histograms of the group mean velocities 
$u=\frac{1}{N}\sum_{i=1}^N v_i$
recorded in $10^5$ time steps of length $10^{-2}$, with $N=5$ (left panel), $N=7$ (middle panel) and $N=12$ individuals (right panel).
The system does not prefer any particular state for $N=5$.
Two quasi-stable states of ordered collective motion are easily recognizable for $N=12$.}}
\label{fig:example1}
\end{figure}

\begin{figure}
\begin{center}
\epsfig{figure=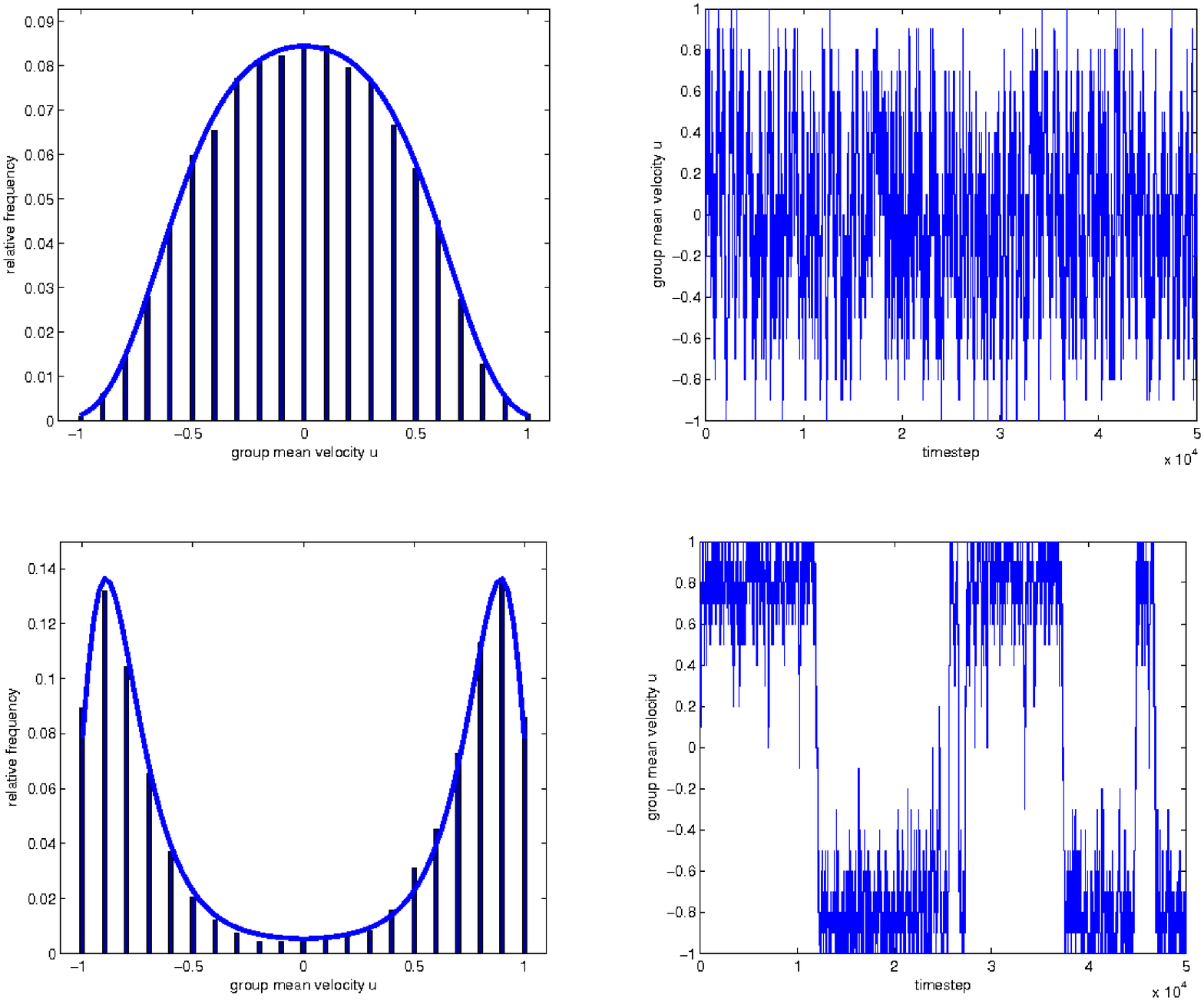, width=1.0\linewidth}
\end{center}
\caption{{\it An example of the behaviour of the model~$\fref{model0}$--$\fref{gammaLR}$ with global interactions ($\w\equiv 1$):
the large noise case (top) with $N=20$, $b=1$ and $\gamma_0=1.3$
and small noise case (bottom) with $N=20$, $b=1$ and $\gamma_0=0.3$.
Left are the histograms of the group mean velocities $u$ recorded in $10^5$ time steps of length $10^{-2}$,
compared to the plot of the (properly scaled) stationary solution $p_s$
of the corresponding Fokker-Planck equation (solid line).
Right are the plots of the temporal evolution of the group mean velocity $u$
during $5\times 10^4$ timesteps.
In the small noise case, one can clearly distinguish the two quasi-stationary states
and observe the switching between them.}}
\label{fig:example2}
\end{figure}

With $\xi(s) = s^2$, the turning rate ``from the right to the left'', $\gamma^{R\to L}$,
and the rate for the opposite turn, $\gamma^{L\to R}$, are given by
\begin{eqnarray}
   \gamma^{R\to L}_i &=& \gamma_0 + b \, (1-u_i^{loc})^2
   \quad\mbox{for the switch from } v_i = 1 \mbox{ to } v_i = -1 \,,
   \label{gammaRL} \\
   \gamma^{L\to R}_i &=& \gamma_0 + b \, (1+u_i^{loc})^2
   \quad\mbox{for the switch from } v_i = -1 \mbox{ to } v_i = 1 \,.
   \label{gammaLR}
\end{eqnarray}
This velocity jump process describes the tendency of the individuals
to align their velocities to the average velocity of their neighbors.
Despite its relative simplicity, the model provides similar predictions as
the Vicsek and Czir\'ok model~\cite{Czirok:1999:CMS} and its modification~\cite{Yates:2009:INC}
and is in qualitative agreement with the experimental observations made in~\cite{Buhl:2006:DOM}:
the transition to ordered motion as $N$ grows (Figure~\ref{fig:example1})
and the density-dependent switching behaviour between the ordered states
(Figure~\ref{fig:example2}, bottom).

%

\section{Analysis of the individual based model with global interactions}\label{sec:Fokker-Planck}

In this section, we simplify the individual-based model by assuming 
$\w\equiv 1$ in (\ref{u_i^loc}), i.e.
$u_i^{loc} = u$ for all $i=1,\dots,N$, where
\begin{equation}
    u(t) := \frac{1}{N} \sum_{i=1}^N v_i(t)  
    = \frac{2 r(t) - N}{N}\,,
    \label{uglobal}
\end{equation}
where $r(t)$ is the number of individuals which are going to the right
at time $t$. 
Using this simplification, we will obtain an explicit formula
for the mean switching time between the ordered states.
However, for the derivation of the kinetic decription and its analysis
(Section~\ref{sec:Kinetic}),
we will allow general weights $\w$, imposing only the assumption
\ref{A1} and a slightly reinforced version of \ref{A2}.

Let $p(r,t)$ be the probability that $r$ individuals 
move to the right (i.e., with velocity $1$) at time $t\geq 0$. 
It satisfies the master equation
\begin{eqnarray}
   \part{}{t} {p(r,t)} &=& (r+1) \, \gamma^{R\to L} p(r+1,t) 
                    - r \, \gamma^{R\to L} p(r,t) 
   \nonumber \\
		&+& (N-r+1) \, \gamma^{L\to R} p(r-1,t) 
		- (N-r) \, \gamma^{L\to R} p(r,t) \,,
   \label{masterglobal}
\end{eqnarray}
where 
$\gamma^{R\to L}$ and $\gamma^{L\to R}$ are given by
(\ref{gammaRL}) and (\ref{gammaLR}), respectively. The subscript
$i$ in (\ref{gammaRL})--(\ref{gammaLR}) is dropped in (\ref{masterglobal}) 
because all individuals have the same turning rates. Using the system size 
expansion \cite{vanKampen:2007:SPP} and the definition (\ref{uglobal})
of the average velocity $u(t)$, we obtain the
following Fokker-Planck equation 
\begin{equation}
   \label{Fokker-Planck}
   \part{p(u,t)}{t} 
   = \part{}{u}
   \Big(
    2 u \left[ \gamma_0 - b (1-u^2) \right]
   p(u,t) \Big) 
   + \partk{}{u}{2}
   \left(
   \frac{2}{N} \left[\gamma_0 + b(1-u^2)\right] 
   p(u,t) \right) \, ,
\end{equation}
where $p(u,t)$ is the probability distribution function of the average velocity
(\ref{uglobal}) at time $t$. The stationary solution $p_s$ of~\fref{Fokker-Planck} 
is 
\(   \label{p_s}
   p_s(u) = C \exp[-\Phi_N(u)]
\)
where $C$ is the normalization constant and the potential $\Phi_N$ is given by
\begin{equation}
   \Phi_N(u)= -\frac{N}{2}u^2 + 
   \left( 1 - \frac{\gamma_0 N}{b} \right)
   \ln\left(\gamma_0 + b(1-u^2)\right) \,.
\label{potentialPhiN}
\end{equation}
The comparision of the stationary probability distribution function 
$p_s$ with the results obtained by long-time simulation of the 
stochastic individual-based model is shown in Figure~\ref{fig:example2}.
Differentiating (\ref{potentialPhiN}) twice, we obtain
\begin{equation}
   \Phi_N''(0)
   = 
   \frac{N}{(\gamma_0 + b)^2}
   \left( 
         \left( \gamma_0 - \frac{b}{N} \right)^2
         - 
	 \left( b + \frac{b}{N} \right)^2
   \right)
\label{derphizero}
\end{equation}
Consequently, we distinguish the following two cases:
\begin{enumerate}[leftmargin = 1 cm, label = {(}\arabic*{)} ]
\item
{\it Large noise}: If 
$\;\displaystyle \frac{\gamma_0}{b} \geq 1 + \frac{2}{N} \;$, then
$\Phi_N$ has the global minimum at $u=0$.
\item
{\it Small noise}: If 
$\;\displaystyle \frac{\gamma_0}{b} < 1 + \frac{2}{N} \;$, then
$\Phi_N$ has a local maximum at $u=0$. The only local and 
global minima are at $\rule{0pt}{5mm} \pm u_s$ where
\[
   u_s = \sqrt{1 + \frac{2}{N} - \frac{\gamma_0}{b}}.
\]
\end{enumerate}
It is worth noting that $u_s > 1$ if $N$ is small,
namely for $N < 2 b/\gamma_0.$ This is a consequence
of approximations made during the derivation of the
Fokker-Planck equation (\ref{Fokker-Planck}).
  
In the large noise case, the system prefers the disordered state $u=0$,
while, in the small noise case, the system has two preferred ordered 
states $\pm u_s$ where it spends most of the time and eventually 
switches between them (see Figure~\ref{fig:example2}, bottom panel).

Using Kramers theory  \cite{Hanggi:1990:RTF,Gillespie:1992:MPI}, the mean 
switching time $\tau_N$ between the states $-u_s$ and $u_s$ can be approximated
as
\( \label{mean_sw_time}
   \tau_N(-u_s \mapsto u_s) 
   \approx \frac{N \pi}{\gamma_0 + b} 
   \frac{\exp(\Phi_N(0)-\Phi_N(-u_s))}{\sqrt{-\Phi_N''(0)\Phi_N''(-u_s)}}  \,,
\)
which has an exponential asymptotic with respect to large $N$ given by
\[
    \tau_N(-u_s \mapsto u_s) \approx \frac{2 \pi}{b-\gamma_0} \sqrt{\frac{\gamma_0}{2(\gamma_0+b)}}
    \exp\left\{ N\left[ \frac{b-\gamma_0}{2} - \frac{\gamma_0}{b}\ln\left(\frac{2\gamma_0}{\gamma_0+b}\right)\right]\right\} \,.
\]
This is in agreement with the experimental observations,~\cite{Buhl:2006:DOM},
as well as with the modified Czirok-Vicsek model of~\cite{Yates:2009:INC},
where the mean switching time is as well exponential in $N$.
Finally, it is interesting to note that with the transform $b = \gamma_0\bar{b}$
in~\fref{mean_sw_time}, one has
\[
   \tau(N,\gamma_0,b) = \gamma_0\bar{\tau}(N,\bar{b}) \,,
\]
i.e., if $b/\gamma_0$ is kept fixed, the mean switching time scales linearly with $\gamma_0$.

In the numerical experiment shown in Figure~\ref{fig:example2} bottom
(small noise case with $\gamma_0 = 0.3$, $b=1$ and $N=20$ agents),
we have two metastable states located approximately at
$u_s = \pm 0.894$. The estimate mean turning time given by formula \fref{mean_sw_time}
is $\tau_N(-u_s \mapsto u_s) = 61.1$.
Performing $10^6$ time steps of length $10^{-2}$, the observed mean switching time
(defined as a mean transition time between the states $v=-0.8$ and $v=0.8$
or vice versa) was $58.8$, showing a very good agreement.

\def\G{\mathcal{G}}
\def\L{\mathcal{L}}
\section{Kinetic description}\label{sec:Kinetic}
In this section we derive the kinetic description of the system of $N$
interacting agents and formally pass to the limit $N\to\infty$
to obtain the corresponding kinetic equation.
For this, we have to accept a slight reinforcement of the assumption~\ref{A2} on $\w$,
namely,
\begin{enumerate}[leftmargin = 1.6 cm, label={[{A2'}]}]
\item $\w > 0$ on the interval $[0,r)$ for some $r>0$. \label{A2'}
\end{enumerate}

\noindent
The state of the system of $N$ agents at time $t\geq 0$ is described by the probability density
function $p^N(t, x_1,v_1,\dots,x_N,v_N)$ of finding the $i$-th agent in position $x_i\in\Omega$
with velocity $v_i\in\{-1,1\}$, for $i=1,\dots,N$.
When convenient, we will use the abbreviation $p^N = p^N(t,\boldx,\boldv)$,
with $\boldx=(x_1,\dots,x_N)\in\Omega^N$ and 
$\boldv=(v_1,\dots,v_N)\in\setv{N}$.
The probability density $p^N(t,\boldx,\boldv)$ evolves according to
\(  \label{kin1}
   \part{}{t}p^N + \sum_{i=1}^N v_i\part{}{x_i} p^N = -\L[p^N] + \G[p^N] \,,
\)
where $\L[p^N]$ (resp. $\G[p^N]$) is the loss (resp. gain) term
corresponding to the velocity jumps. Using (\ref{gammaRL})--(\ref{gammaLR}),
the rate of the switch $v_i \mapsto -v_i$ is for $i\in\{1,\dots,N\}$ 
given by
\[
   \gamma_i(\boldx,\boldv) 
   = 
   \gamma_0 + b\,\big(v_i - u[x_i; \boldx,\boldv] \big)^2 \,,
\]
where $u[x_i; \boldx,\boldv] = u_i^{loc}$ is defined by~\fref{u_i^loc},
i.e. $u[x_i; \boldx,\boldv]$
is the local average velocity seen by an agent 
located at $x_i \in \Omega$, based on the system configuration 
$[\boldx,\boldv]$. Consequently, the loss term is given by
\(  \label{loss}
   \L[p^N](t,\boldx,\boldv) = \sum_{i=1}^N 
   \gamma_i(\boldx,\boldv) \, p^N(t,\boldx,\boldv) \,.
\)
We define the operator $M_i: \setv{k} \to \setv{k}$ for
$k \in \{1,\dots,N\}$ and $i \in \{1,\dots,k\}$ by
\begin{equation}
M_i(\boldv) = (v_1,\dots,v_{i-1},-v_i,v_{i-1},\dots,v_k),
\label{defMi}
\end{equation}
i.e. $M_i(\boldv)$ denotes the velocity vector created
from $\boldv$ by changing the sign of its $i$-th component.
Then the gain term $\G[p^N](t,\boldx,\boldv)$ is given by
\(   \label{gain}
   \G[p^N](t,\boldx,\boldv) 
   = 
   \sum_{i=1}^N 
   \gamma_i(\boldx,M_i(\boldv)) \, p^N(t,\boldx,M_i(\boldv)) \, ,
\)
where
\[
   \gamma_i(\boldx,M_i(\boldv)) 
   = \gamma_0 
    + b
      \,
      \big(\!-v_i - u[x_i; \boldx, M_i(\boldv)] \big)^2 \,.
\]
Consequently, the right hand side of~\fref{kin1} is 
\(   \label{coll1}
   (\G - \L)\left[p^N\right](\boldx,\boldv) 
   &=& \gamma_0 \left( \sum_{i=1}^N 
   p^N(\boldx,M_i(\boldv)) - Np^N(\boldx,\boldv)\right) \\
    &+& b \sum_{i=1}^N \left(
        \,
      \big(\!-v_i - u[x_i; \boldx,M_i(\boldv)] \big)^2
      \,  p^N(\boldx,M_i(\boldv))
      - \big(v_i - u[x_i; \boldx,\boldv] \big)^2
      \, p^N(\boldx,\boldv) \right),
   \nonumber
\)
where we dropped the dependance on time $t$ to simplify the notation.
Finally, we postulate the so-called \emph{indistinguishability-of-particles}:
We only consider solutions $p^N$ that are indifferent to permutations
of their $(x_i,v_i)$-arguments. Such solutions are admissible,
since the equation~\fref{kin1}
with the collision operator~\fref{coll1} is as well indifferent
with respect to interchange of the $(x_i,v_i)$-pairs.

\subsection{Derivation of the BBGKY hierarchy}
To derive an analogue of what is called the \emph{BBGKY hierarchy}
in the classical kinetic theory of gases (see, for instance, \cite{Cercignani:1994:MTD}),
we define, for $k=1,\dots,N$, the $k$-agent marginals
\( \label{marginal}
   p^{N,k}(\boldx^k,\boldv^k) 
   = \sum_{\bar \boldv\in\setv{N-k}} \int_{\Omega^{N-k}}
         p^N(\boldx^k,\bar \boldx,\boldv^k,\bar \boldv) \d\bar \boldx \,,
    \qquad \boldx^k\in\Omega^k \,,\; \boldv^k\in\setv{k}\,.
\)
In what follows, the coordinates of $\boldx^k$ and $\bar\boldx$ will be denoted as
\[
    \boldx^k = (x^k_1, x^k_2, \dots, x^k_k) \,,\qquad\qquad\qquad
    \bar\boldx = (\bar x_1, \bar x_2, \dots, \bar x_{N-k}) \,,
\]
that is, $\boldx = (\boldx^k, \bar\boldx) = (x^k_1, x^k_2, \dots, x^k_k, \bar x_1, \bar x_2, \dots, \bar x_{N-k})$.
The same notational convention is used for velocities, i.e.,
$\boldv = (\boldv^k, \bar\boldv) = (v^k_1, v^k_2, \dots, v^k_k, \bar v_1, \bar v_2, \dots, \bar v_{N-k})$.
Note that, due to the indistinguishability-of-particles, the marginals are well defined
and indifferent with respect to permutations of the pairs of arguments $(x_i,v_i)$,
$i=1,\dots,k$. Integrating \fref{kin1}, we obtain
\begin{equation}
    \part{}{t} p^{N,k}(\boldx^k,\boldv^k) 
    + \sum_{i=1}^k v_i \part{}{x_i} p^{N,k}(\boldx^k,\boldv^k) =
      \sum_{\bar \boldv\in\setv{N-k}} 
      \int_{\Omega^{N-k}} 
      (\G-\L)[p^N](\boldx^k,\bar \boldx,\boldv^k,\bar \boldv) \d\bar \boldx\,.
\label{evolpNk}
\end{equation}
Substituting (\ref{coll1}) into the right-hand side of (\ref{evolpNk}), 
we get
\begin{eqnarray}
    \sum_{\bar \boldv\in\setv{N-k}} \int_{\Omega^{N-k}} 
    (\G-\L)[p^N]\!\!\!\!\!\!\!\!\!\!\!\!\!\!\!\!
    &&(\boldx^k,\bar \boldx,\boldv^k,\bar \boldv) \d\bar \boldx \nonumber \\
    &=& 
    \gamma_0 \, \sum_{i=1}^k 
     (p^{N,k}(\boldx^k,M_i(\boldv^k)) - p^{N,k}(\boldx^k,\boldv^k))
\nonumber \\
&+& b \sum_{\bar \boldv\in\setv{N-k}} \int_{\Omega^{N-k}} 
    \sum_{i=1}^k 
     \Bigl\{ 
      \big(\!-v_i^k - 
      u[x_i; \boldx^k,\bar \boldx,M_i(\boldv^k),\bar \boldv] \big)^2 
             p^N(\boldx^k,\bar \boldx,M_i(\boldv^k),\bar \boldv) \Bigr.\; 
	     \nonumber \\
&& \qquad\qquad	\qquad\qquad \qquad     
     \Bigl. - \,
     \big(v_i^k - u[x_i; \boldx^k,\bar \boldx,\boldv^k,\bar \boldv] \big)^2
             p^N(\boldx^k,\bar \boldx, \boldv^k,\bar \boldv)
             \Bigr\}\, \nonumber \\
&=& (\gamma_0+b)\, \sum_{i=1}^k 
    (p^{N,k}(\boldx^k,M_i(\boldv^k)) - p^{N,k}(\boldx^k,\boldv^k)) 
   \nonumber \\ 
&+&  
   2b \!
   \sum_{\bar \boldv\in\setv{N-k}} 
   \int_{\Omega^{N-k}} \sum_{i=1}^k
      v^k_i
      \left( 
      u[x_i; \boldx^k,\bar \boldx,M_i(\boldv^k),\bar \boldv]
      \,  p^N(\boldx^k,\bar \boldx,M_i(\boldv^k),\bar \boldv)
       \right.  \nonumber\\
	 && \qquad\qquad	\qquad\qquad \qquad 
	 +\left.
	 u[x_i; \boldx^k,\bar \boldx,\boldv^k,\bar \boldv] 
	 \,
	 p^N(\boldx^k,\bar \boldx,\boldv^k,\bar \boldv) \right)
	 \d\bar \boldx
\nonumber \\
&+&  
 b \!
   \sum_{\bar \boldv\in\setv{N-k}} 
   \int_{\Omega^{N-k}} \sum_{i=1}^k 
       \left( 
       \big( u[x_i; \boldx^k,\bar \boldx,M_i(\boldv^k),\bar \boldv] \big)^2
       \, p^N(\boldx^k,\bar \boldx,M_i(\boldv^k),\bar \boldv)
       \right. \nonumber \\
     && \qquad\qquad	\qquad\qquad \qquad -
	 \left.
       \big( u[x_i; \boldx^k,\bar \boldx,\boldv^k,\bar \boldv] \big)^2     
     \, p^N(\boldx^k,\bar \boldx,\boldv^k,\bar \boldv) \right)
	 \d\bar \boldx.
	 \label{derBBGKY}
\end{eqnarray}
Since we are interested in the limit $N\to\infty$ with $k$ fixed,
we can rewrite the definition (\ref{u_i^loc}) as follows
\begin{eqnarray}
   u[x_i; \boldx^k,\bar \boldx,\boldv^k,\bar \boldv] &=& 
       \frac{\sum_{m=1}^k \w(|x^k_m-x_i|) v^k_m
       + \sum_{m=1}^{N-k} \w(|\bar x_m - x_i|) \bar v_m}
            {\sum_{m=1}^k \w(|x^k_m-x_i|) 
	    + \sum_{m=1}^{N-k} \w(|\bar x_m - x_i|)} \nonumber\\
      &=& 
	  u[x_i;\bar \boldx,\bar \boldv] + \O\left(\frac{k}{N}\right),
	  \,
\label{averu}	  
\end{eqnarray}
where we define
\[
    u[z;\bar \boldx,\bar \boldv]
    =
    \frac{\sum_{m=1}^{N-k} 
    \w(|\bar x_m - z|) \bar v_m}{\sum_{m=1}^{N-k} \w(|\bar x_m-z|)}  \,
    \qquad  \mbox{for} \; z \in \Omega.  
\]     
Substituting (\ref{averu}) in (\ref{derBBGKY}) and (\ref{evolpNk}), 
we obtain the \emph{BBGKY hierarchy}
\(
    \!\!\!\!\part{}{t} p^{N,k}(\boldx^k,\boldv^k) 
    + \sum_{i=1}^k v_i \part{}{x_i} p^{N,k}(\boldx^k,\boldv^k) 
     \!\!\!&=&\!\!\!
     (\gamma_0+b)\, \sum_{i=1}^k 
    \left(p^{N,k}\!\left(\boldx^k,M_i\big(\boldv^k\big)\right) 
    - p^{N,k}\!\left(\boldx^k,\boldv^k\right)\right) 
     \nonumber\\
      &+&\!\!\! 
      2b \sum_{i=1}^k v^k_i 
     \left(q^{N,k}\!\left(x_i; \, \boldx^k,M_i\big(\boldv^k\big)\right) 
      + q^{N,k}\!\left(x_i; \, \boldx^k,\boldv^k \right)\right)
      \label{BBGKY} \\
      &+&\!\!\!  b \sum_{i=1}^k 
      \left(r^{N,k}\!\left(x_i; \, \boldx^k,M_i\big(\boldv^k\big)\right) 
      - r^{N,k}\!\left(x_i; \, \boldx^k,\boldv^k\right)\right) 
      + \O\left(\frac{k}{N}\right) \nonumber ,
\)
where
\(   \label{q}
    q^{N,k}\!\left(z; \, \boldx^k,\boldv^k\right) 
    = \sum_{\bar \boldv\in\setv{N-k}} 
    \int_{\Omega^{N-k}}
    u[z; \bar \boldx,\bar \boldv] 
    \, p^{N}\big(\boldx^k,\bar \boldx,\boldv^k,\bar \boldv\big) 
    \d\bar \boldx \,.
\)
and
\(  \label{r}
    r^{N,k}\!\left(z; \, \boldx^k,\boldv^k\right) 
    = \sum_{\bar v\in\setv{N-k}} 
    \int_{\Omega^{N-k}}
     u[z; \bar \boldx,\bar \boldv]^2 
     \, p^{N}\big(\boldx^k,\bar \boldx,\boldv^k,\bar \boldv\big) 
     \d\bar \boldx \,. 
\)

\subsection{Passage to the limit $N\to\infty$}\label{subsect:Ntoinfty}
The usual procedure of deriving the mean field equation
is to write the \emph{BBGKY hierarchy}~\fref{BBGKY} in terms of $p^{N,k}$
and pass to the limit $N\to\infty$ to obtain the so-called 
\emph{Boltzmann hierarchy} for $p^k := \lim_{N\to\infty} p^{N,k}$ \cite{Sznitman:1991:TPC}.
Then, one shows that the Boltzmann hierarchy admits solutions
generated by the \emph{molecular chaos ansatz} (see below).
In our case, however, this strategy cannot be pursued;
although we could derive uniform estimates allowing us to
pass to the limit $N\to\infty$ in the BBGKY hierarchy,
we are not able to express the limiting marginals
$q^k := \lim_{N\to\infty} q^{N,k}$ and $r^k := \lim_{N\to\infty} r^{N,k}$
in terms of $p^k$. Consequently, we have no clue what the 
correct molecular chaos ansatz for $q^k$ and $r^k$ should be.

Instead, we assume the propagation of chaos already at the level of the BBGKY hierarchy,
before passing to the limit $N\to\infty$: we assume that, for large $N$,
$p^N$ is well approximated by the product of the limiting one-particle marginals
$p := \lim_{N \to \infty} p^{N,1}$, i.e.
\( \label{molecular chaos}
    p^N(t,\boldx,\boldv) \approx \prod_{i=1}^N p(t,x_i,v_i)
         \qquad\mbox{for all }   
	 t\geq 0\,,\;\boldx\in\Omega^N\,,\; \boldv\in\setv{N} \,.
\)
This corresponds to vanishing statistical dependence (correlations) between the agents
as $N\to\infty$ and is the usual phenomenon observed in systems of interacting particles,
see for instance~\cite{Cercignani:1994:MTD} in the context of classical kinetic theory
or~\cite{Haskovec:2009:SPA, Haskovec:2010:CAS} in the context of biological systems.
Moreover, if one interprets $q^{N,k}$ (resp. $r^{N,k}$) as the first (resp. second) order
moment of $p^N$ with respect to $u[z;\bar\boldx,\bar\boldv]$, then one can understand~\fref{molecular chaos}
as the moment closure assumption for the non-closed system of moments generated by~\fref{kin1}.

The essential point is that now we may insert~\fref{molecular chaos}
into~\fref{q} and~\fref{r} to obtain explicit expressions for $q^k$ and $r^k$
in terms of $p$:
\(  \label{lim-q}
    && q^k(z; \boldx^k,\boldv^k) = \lim_{N\to\infty} q^{N,k}\!\left(z; \, \boldx^k,\boldv^k\right) = \\
       &&   \left(\prod_{i=1}^k p(x^k_i,v^k_i) \right) \times
        \lim_{N\to\infty} \sum_{\bar\boldv\in\setv{N-k}} \int_{\Omega^{N-k}}
          \frac{\sum_{m=1}^{N-k} \w(|\bar x_m - z|) \bar v_m}{\sum_{m=1}^{N-k} \w(|\bar x_m - z|)}
          \prod_{i=1}^{N-k} p(\bar x_i,\bar v_i) \d\bar\boldx \,,  \nonumber
\)
and
\(  \label{lim-r}
    && r^k(z;\boldx^k,\boldv^k) = \lim_{N\to\infty} r^{N,k}\!\left(z; \, \boldx^k,\boldv^k\right) = \\
       &&   \left(\prod_{i=1}^k p(x^k_i,v^k_i) \right) \times
       \lim_{N\to\infty} \sum_{\bar\boldv\in\setv{N-k}} \int_{\Omega^{N-k}}
          \left( \frac{\sum_{m=1}^{N-k} \w(|\bar x_m - z|) \bar v_m}{\sum_{m=1}^{N-k} \w(|\bar x_m - z|)}\right)^{\!\!2}\,
          \prod_{i=1}^{N-k} p(\bar x_i,\bar v_i) \d\bar\boldx \,. \nonumber
\)

\subsubsection*{Study of the limit $N\to\infty$ in~\fref{lim-q} and~\fref{lim-r}}
We start by setting $k=1$, which is the case considered in Section~\ref{subsect:Derivation}.
To simplify the notation, we drop the bars over $x$ and $v$, and, without loss of generality, choose $z=0$.
First, we explore the symmetry of the expression~\fref{lim-q} as follows:
\[
   \sum_{\boldv\in\setv{N}} \int_{\Omega^{N-1}} \!\!\! && \!\!\!\!\!\!\!\!\!
          \frac{\sum_{m=1}^{N-1} \w(x_m) v_m}{(N-1) S_{N}(\boldx)}
          \prod_{i=1}^{N-1} p(x_i,v_i) \d\boldx \\
  &=& \sum_{m=1}^{N-1} \int_{\Omega^{N-1}} \frac{\w(x_m)}{(N-1) S_{N}(\boldx)}
         \left[ p(x_m,1) - p(x_m,-1) \right]
     \sum_{\boldv\in\setv{N-2}} \prod_{i\neq m} p(x_i,v_i) \d\boldx \\
  &=& \sum_{m=1}^{N-1} \int_{\Omega^{N-1}} \frac{\w(x_m)}{(N-1) S_{N}(\boldx)}
         j(x_m) \prod_{i\neq m} \varrho(x_i) \d\boldx \\
  &=& \int_{\Omega^{N-1}} \frac{\w(x_1)}{S_{N}(\boldx)}
      j(x_1) \prod_{i=2}^{N-1} \varrho(x_i) \d\boldx \,,
\]
with the notation
\[
   \varrho(x) &:=& p(x,1) + p(x,-1) \,,\\
   j(x) &:=& p(x,1) - p(x,-1) \,,
\]
and
\(   \label{S_N}
   S_{N}(\boldx) := \frac{1}{N-1} \sum_{i=1}^{N-1} w(x_i) \,.
\)
Due to the normalization of $p^N$, we have $\int_\Omega \varrho(x) \d x = 1$.
Consequently, in what follows we denote by $P_\varrho(t)$ the time dependent probability measure
corresponding to the probability density $\varrho(t)$.
Since $\w$ is bounded and nonnegative by assumption~\ref{A1}, it is integrable with respect to~$P_\varrho$ and we may define
\(   \label{I}
   I:= \int_\Omega w(x) \d P_\varrho(x) \geq 0\,.
\)
The forthcoming analysis will be performed given the assumption $I>0$;
the case $I=0$ will be discussed in Remark~\ref{Remark:I}.

\begin{lemma} \label{lemma:lim-q}
Let $P_\varrho$ be a probability measure on $\Omega$ with density $\varrho$
and $j\in L^1(\Omega)$ such that $|j| \leq \rho$.
Let $\w:\Omega\to [0,\infty)$ with $w\in L^\infty(\Omega)$ be such that the integral $I$ defined by~$\fref{I}$ is positive.
Define
\[
   Q_N := \int_{\Omega^{N-1}} \frac{\w(x_1)j(x_1)}{S_N(\boldx)} \d x_1 \prod_{i=2}^{N-1}\d P_\varrho(x_i) \,.
\]
Then
\[
   \lim_{N\to\infty} Q_N = \frac{1}{I} \int_\Omega \w(y) j(y)  \d y \,.
\] 
\end{lemma}

\begproof
We can treat $\w(y)$ as a random variable with respect to the probability measure $P_\varrho(y)$.
The essential tool of the proof is the law of large numbers,
which states that $S_N(\boldx)$ converges to $I$ in measure,
in the sense that for each $\eps>0$,
\(   \label{law_large_numbers}
   \lim_{N\to\infty} P_\varrho^{N-1}\left(\{\boldx\in\Omega^{N-1};\; |S_N(\boldx)-I| > \eps\}\right) = 0 \,,
\)
where $P^{N-1}_\varrho$ denotes the $(N\!-\!1)$-fold tensor product of the probability measures $P_\varrho$.
Moreover, the existence of the $m$-th order moment of $\w$ with respect to $P_\varrho$,
\[
    \int_\Omega |\w(y)|^m \d P_\varrho(y) < \infty \,,
\]
implies the rate of convergence (see~\cite{Baum:1965:CRL})
\(   \label{rate of convergence}
   \lim_{N\to\infty} (N-1)^{m-1} P_\varrho^{N-1}\left(\{\boldx\in\Omega^{N-1};\; |S_N(\boldx)-I| > \eps\}\right) = 0 \,.
\)
Let us denote by $A_N(\eps)$ the set $\{\boldx\in\Omega^{N-1},\; |S_N(\boldx)-I| < \eps\}$
and by $A_N^c(\eps)$ its complement in $\Omega^{N-1}$.
Choosing $0 < \eps < I/2$ and $\boldx\in A_N(\eps)$, we have the estimate
\[
    \left| \frac{w(x_1)}{S_N(\boldx)} - \frac{w(x_1)}{I} \right| \leq
      \frac{2\eps}{I^2} \w(x_1) \,.
\]
Consequently,
\[
   \int_{A_N(\eps)} \left| \frac{w(x_1)}{S_N(\boldx)} - \frac{w(x_1)}{I} \right| \d P_\varrho^{N-1}(\boldx)
   \leq \frac{2\eps}{I^2} \int_{\Omega^{N-1}} \w(x_1) \d P_\varrho^{N-1}(\boldx) = \frac{2\eps}{I}  \,.
\]
On the other hand, the integral over $A_N^c(\eps)$ is estimated with
\[
   \int_{A_N^c(\eps)} \left| \frac{w(x_1)}{S_N(\boldx)} - \frac{w(x_1)}{I} \right| \d P_\varrho^{N-1}(\boldx)
   \leq
   \int_{A_N^c(\eps)} \frac{w(x_1)}{I} \d P_\varrho^{N-1}(\boldx) +
   \int_{A_N^c(\eps)} \frac{w(x_1)}{S_N(\boldx)} \d P_\varrho^{N-1}(\boldx)  \,.
\]
The first term on the right hand side converges to zero as $N\to\infty$ by~\fref{law_large_numbers} and the boundedness of $\w$.
Using~\fref{S_N}, the second term is estimated by
\[
   \int_{A_N^c(\eps)} \frac{w(x_1)}{S_N(\boldx)} \d P_\varrho^{N-1}(\boldx) \leq (N-1) P_\varrho^{N-1}(A_N^c(\eps))
\]
and vanishes as $N\to\infty$ due to~\fref{rate of convergence} with $m=2$.
Consequently, we have shown that the term
\[
    \left|Q_N - \frac{1}{I}\int_\Omega \w(y) j(y) \d y\right| &=&
    \left| \int_{\Omega^{N-1}} \left( \frac{w(x_1)}{S_N(\boldx)} - \frac{w(x_1)}{I} \right) j(x_1) \d x_1 \prod_{i=2}^{N-1} \d P_\varrho(x_i)\right| \\
    &\leq& \int_{\Omega^{N-1}} \left| \frac{w(x_1)}{S_N(\boldx)} - \frac{w(x_1)}{I} \right| \d P_\varrho^{N-1}(\boldx)
\]
can be made arbitrarily small, for sufficiently large $N$.
\endproof

\noindent
The formula~\fref{lim-r} can be analysed in a similar way as we did for~\fref{lim-q}.
Namely, we exploit its symmetry as follows,
\(
   && \sum_{\boldv\in\setv{N-1}} \int_{\Omega^{N-1}}
          \left( \frac{\sum_{m=1}^{N-1} \w(x_m) v_m}{(N-1) S_N(\boldx)} \right)^{\!\!2} \,
          \prod_{i=1}^{N-1} p(x_i,v_i) \d\boldx \nonumber\\
   &=& \frac{1}{(N-1)^2} \sum_{i=1}^{N-1} \sum_{m\neq i} \int_{\Omega^{N-1}}
         \frac{\w(x_i)\w(x_m)j(x_i) j(x_m)}{S_N^2(\boldx)} \d x_i \d x_m \prod_{k\neq i,m} \d P_\varrho(x_k) \nonumber\\
     &&  + \frac{1}{(N-1)^2} \sum_{m=1}^{N-1} \int_{\Omega^{N-1}}
          \frac{\w(x_m)^2}{S_N^2(\boldx)} \d P_\varrho^{N-1}(\boldx) \nonumber\\
   &=&    \frac{N-2}{N-1} \int_{\Omega^{N-1}}
         \frac{\w(x_1)\w(x_2)j(x_1) j(x_2)}{S_N^2(\boldx)} \d x_1 \d x_2 \prod_{k=3}^{N-1} \d P_\varrho(x_k) \nonumber\\
     &&  + \frac{1}{N-1} \int_{\Omega^{N-1}}
          \frac{\w(x_1)^2}{S_N^2(\boldx)} \d P_\varrho^{N-1}(\boldx) \,. \label{term-r}
\)
The limiting behaviour of the first term is studied in the following Lemma:

\begin{lemma} \label{lemma:lim-r}
With the assumptions and notation of Lemma~\ref{lemma:lim-q},
and defining
\[
   R_N := \int_{\Omega^{N-1}} \frac{\w(x_1)\w(x_2)j(x_1) j(x_2)}{S_N^2(\boldx)} \d x_1 \d x_2 \prod_{k=3}^{N-1} \d P_\varrho(x_k) \,,
\]
we have
\[
   \lim_{N\to\infty} R_N = \left( \frac{1}{I} \int_\Omega w(y) j(y) \d y \right)^2 \,.
\] 
\end{lemma}

\begproof
We follow the lines of the proof of Lemma~\ref{lemma:lim-q}.
Defining again $A_N(\eps) :=\{\boldx\in\Omega^{N-1},\; |S_N(\boldx)-I| < \eps\}$
and $A_N^c(\eps):=\Omega^{N-1} \setminus A_N(\eps)$, with $0< \eps < I/2$, we derive the estimates
\[
    \int_{A_N(\eps)} \left| \frac{\w(x_1)\w(x_2)}{S_N^2(\boldx)} - \frac{\w(x_1)\w(x_2)}{I^2} \right| \d P_\varrho^{N-1}(\boldx)
      \leq \frac{10\eps}{I}
\]
and
\[
   \int_{A_N^c(\eps)} \left| \frac{\w(x_1)\w(x_2)}{S_{N-1}^2(\boldx)} - \frac{\w(x_1)\w(x_2)}{I^2} \right| \d P_\varrho^{N-1}(\boldx)
   &\leq&
   \int_{A_N^c(\eps)} \frac{\w(x_1)\w(x_2)}{I^2} \d P_\varrho^{N-1}(\boldx) \\
     &+& (N-1)^2 \left[P_\varrho^{N-1}(A_N^c(\eps)\right]^2 \,,
\]
where the first term vanishes in the limit $N\to\infty$ due to~\fref{law_large_numbers}
and the second one by~\fref{rate of convergence} with $m=3$.
\endproof

\noindent
By a slight modification of the above Lemma, one obtains the limit of the second term of~\fref{term-r}, namely
\[
   \lim_{N\to\infty} \int_{\Omega^{N-1}}  \frac{\w(x_1)^2}{S_N^2(\boldx)} \d P_\varrho^{N-1}(\boldx) = \frac{1}{I} \int_{\Omega} \w(y)^2 \d P_\varrho(y) \,.
\]
Therefore, the limit as $N\to\infty$ of~\fref{term-r} is
\(  \label{lim-r_2}
   \lim_{N\to\infty} \sum_{\boldv\in\setv{N-1}} \int_{\Omega^{N-1}}
          \left( \frac{\sum_{m=1}^{N-1} \w(x_m) v_m}{(N-1) S_N(\boldx)} \right)^{\!\!2} \,
          \prod_{i=1}^{N-1} p(x_i,v_i) \d\boldx =
   \lim_{N\to\infty} R_N = \left( \frac{1}{I} \int_\Omega w(y) j(y) \d y \right)^2 \!\!,
\)
where we used Lemma~\ref{lemma:lim-r}.

\begin{remark}
Lemmas~\ref{lemma:lim-q} and~\ref{lemma:lim-r} were formulated for the case $k=1$, however,
all the calculations can be easily generalized for any (fixed) value of $k$.
\end{remark}

\subsection{Derivation of kinetic and hydrodynamic description}\label{subsect:Derivation}
Let us denote $p^+(t,x) := p(t,x,1)$ and $p^-(t,x) := p(t,x,-1)$.
We define $u(t,x)$, the continuous analogue of the local average velocity~\fref{u_i^loc}, by
\( \label{u_smooth}
   u(t,x) :=
  \frac{\int_\Omega w(|x-z|) (p^+(t,z)-p^-(t,z))\d z}{\int_\Omega w(|x-z|) (p^+(t,z)+p^-(t,z)) \d z}
  \qquad\mbox{for } x\notin\mathcal{S}_0[p^+,p^-](t) \,,
\)
where
\[
   \mathcal{S}_0[p^+,p^-](t) := \left\{ x\in\Omega;\, \int_\Omega w(|x-z|) (p^+(t,z)+p^-(t,z)) \d z = 0 \right\} \,.
\]
We extend the definition of $u(t,x)$ to the whole domain $\Omega$ by setting $u(t,x)=0$ for $x\in \mathcal{S}_0[p^+,p^-](t)$,
see Remark~\ref{Remark:I}.

The kinetic equation for $p$ is obtained by
setting $k=1$ in the BBGKY-hierarchy~\fref{BBGKY}
with the molecular chaos assumption~\fref{molecular chaos}
and passing to the limit $N\to\infty$.
Using Lemma~\ref{lemma:lim-q} and formula~\fref{lim-r_2}, we obtain
\(  \label{q+r}
   q^1(x; x,v) = p(t,x,v) u(t,x) \,,\qquad\mbox{and}\qquad r^1(x; x,v) = p(t,x,v) u(t,x)^2 \,.
\)
Consequently, \fref{BBGKY} reduces to the following system of two equations
\(
    \partial_t{p^+} + \partial_x{p^+} &=& -[\gamma_0 + b(1-u)^2] p^+ + [\gamma_0 + b(1+u)^2] p^-\,,  \label{kinetic1}\\
    \partial_t{p^-} - \partial_x{p^-} &=& -[\gamma_0 + b(1+u)^2] p^- + [\gamma_0 + b(1-u)^2] p^+\,. \label{kinetic2}
\)
Equivalently, defining the mass density $\varrho$ and flux $j$ by
\(   \label{rho-j}
    \varrho(t,x):= p^+(t,x) + p^-(t,x) \,,\qquad
    j(t,x) := p^+(t,x) - p^-(t,x) \,,
\)
the system can be written in the hydrodynamic description as
\(
    \partial_t{\varrho} + \partial_x{j} &=& 0\,,\label{hydrodynamic1}\\
    \partial_t{j} + \partial_x{\varrho} &=& -2 [\gamma_0 + b(1+u^2)] j + 4b\varrho u\,,\label{hydrodynamic2}\\
    \rule{0pt}{32pt} u(t,x) &:=& \left\{ \begin{array}{cl} 
           \displaystyle\frac{\int_\Omega w(|x-z|) j(t,z)\d z}{\int_\Omega w(|x-z|) \varrho(t,z) \d z} \,,& x\notin\mathcal{S}_0[\varrho](t) \,,\\
           \rule{0pt}{15pt}0 \,, & x\in\mathcal{S}_0[\varrho](t) \,, \end{array}\right. \label{hydrodynamic3}
\)
with
\[
   \mathcal{S}_0[\varrho](t) := \left\{ x\in\Omega;\, \int_\Omega w(|x-z|) \varrho(t,z) \d z = 0 \right\} \,.
\]

\begin{remark}\label{Remark:I}
Due to the assumption~$\ref{A2'}$, we have $p^+(t,x)=p^-(t,x)=0$ on $\mathcal{S}_0[p^+,p^-](t)$.
Therefore, both $q^{N,1}\!\left(x; \, x,v\right)$ in~$\fref{lim-q}$ and $r^{N,1}\!\left(x; \, x,v\right)$ in~$\fref{lim-r}$
are equal to zero by definition. Consequently,
\[
   q^1(x; x,v) = \lim_{N\to\infty} q^{N,1}\!\left(x; \, x,v\right) = 0 \,,\qquad\mbox{and}\qquad
   r^1(x; x,v) = \lim_{N\to\infty} r^{N,1}\!\left(x; \, x,v\right) = 0 \,,
\]
and formulas~$\fref{q+r}$ remain valid irrespective of the particular choice of the value of $u(t,x)$.
This justifies our extension of the definition of $u(t,x)$ by setting it equal to zero for $x\in \mathcal{S}_0[p^+,p^-](t)$.
\end{remark}

\section{Existence of solutions to the kinetic system}\label{sec:Existence}
The main goal of this Section is to prove the following Theorem:

\begin{theorem}\label{thm:existence}
Let $\gamma_0 \geq 0$, $b \geq 0$ and $\w$ satisfy the assumptions~$\ref{A1}$ and~$\ref{A2'}$.
Then, for every $T>0$ and every nonnegative initial datum $(p^+_0, p^-_0)\in L^\infty(\Omega)\times L^\infty(\Omega)$
there exists a nonnegative solution to the kinetic formulation~$\fref{kinetic1}$--$\fref{kinetic2}$
in $L^\infty([0,T]\times\Omega)$.
This also establishes solutions $\varrho=p^++p^-$, $j=p^+-p^-$
of the hydrodynamic formulation~$\fref{hydrodynamic1}$--$\fref{hydrodynamic3}$
with the corresponding initial condition.
\end{theorem}

\def\wj{J}
\def\wrho{R}

\noindent\textbf{Proof of Theorem~\ref{thm:existence}:}
The proof is carried out in three steps and is only sketched here, omitting details
where the techniques are standard. For notational convenience, we will work both
with the kinetic and hydrodynamic representation of the system
and treat $(p^+,p^-)$ and $(\varrho,j)$ as synonyms,
related by~\fref{rho-j}.

\medskip\noindent\textbf{Step 1.}
First, we consider a linearized version of~\fref{kinetic1}--\fref{kinetic2},
where we solve for $p^+$ and $p^-$ given a prescribed $u\in L^\infty([0,T]\times\Omega)$ with $|u|\leq 1$.
This constitutes a strictly hyperbolic system with unique mild, nonnegative solution in $C([0,T]; L^\infty(\Omega))$,
constructed by a standard fixed point iteration (see, for instance,~\cite{Evans:1998:PDE}, Section 7.3).
The solution is given by the Duhamel formula
\(
   p^{+}(t,x) &=& p^{+}(0,x-t) + \int_0^t Q^+[p^{+},p^{-}](s,x+(s-t)) \d s \,, \label{Duhamel1}\\
   p^{-}(t,x) &=& p^{-}(0,x+t) + \int_0^t Q^-[p^{+},p^{-}](s,x-(s-t)) \d s \,, \label{Duhamel2}
\)
with $Q^\pm[p^{+},p^{-}] = \mp [\gamma_0 + b(1-u)^2] p^{+} \pm [\gamma_0 + b(1+u)^2] p^{-}$.
Moreover, for any fixed $T>0$, we have apriori boundedness of $p^+$ and $p^-$ in $L^\infty([0,T]\times\Omega)$,
depending only on the initial condition. Indeed,
denoting $h(t) := \sup_{x\in\Omega} p^{+}(t,x) + \sup_{x\in\Omega} p^{-}(t,x)$,
and remembering the uniform boundedness $|u| \leq 1$, we have
\[
    h(t) \leq h(0) + (\gamma_0 + 4b) \int_0^t h(s) \d s \,,
\]
and the apriori boundedness follows from an application of the Gronwall lemma on the time interval $[0,T]$.

\medskip\noindent\textbf{Step 2.}
We consider a regularized version of~\fref{kinetic1}--\fref{kinetic2}
where $u$ is substituted by $u_\eps$, defined by
\( \label{u_eps}
   u_\eps(x) := \frac{\wj(x)}{\eps+\wrho(x)}\qquad \mbox{ with}\quad
   \wj(x):= \int_\Omega w(|x-z|) j(z) \d z \,,\qquad
    \wrho(x):= \int_\Omega w(|x-z|) \varrho(z) \d z \,.
\)
For any fixed $\eps>0$, a solution is found by the Schauder fixed point iteration on the mean velocity $u_\eps$~\cite{Evans:1998:PDE}.
The compactness of the corresponding Schauder operator is provided by the Arzela-Ascoli theorem.
Indeed, let us take a sequence $u^n$ with $\Norm{u^n}_{L^\infty([0,T]\times\Omega)} \leq 1$
and let $p^{+,n}$ and $p^{-,n}$ be the corresponding mild solutions of the kinetic system,
given by the Duhamel formula~\fref{Duhamel1}--\fref{Duhamel2} with $u^n$ in place of $u$,
and let $\varrho^n = p^{+,n}+p^{-,n}$ and $j^n = p^{+,n}-p^{-,n}$.
As explained in Step 1, for any fixed $T>0$ we have $\Norm{p^{+,n}}_{L^\infty([0,T]\times\Omega)}$,
$\Norm{p^{+,n}}_{L^\infty([0,T]\times\Omega)}$ bounded uniformly with respect to $n$.
Defining the function
\[
    \omega(x) := \int_\Omega |\w(z) - \w(|z-x|)| \d z \qquad\mbox{for } x\in\Omega \,,
\]
one has $\lim_{x\to 0} \omega(x) = 0$ (continuity of translation~\cite{Rudin:1991:FA}).
By H\"older inequality,
\[
    |\wj^n(x)-\wj^n(y)| \leq \Norm{j^n}_{L^\infty(\Omega)} \omega(x-y) \,,
\]
with $\wj^n(x):= \int_\Omega w(|x-z|) j^n(z) \d z$.
Moreover, we have the uniform boundedness
\[
    |\wj^n(x)| \leq \Norm{j^n}_{L^\infty(\Omega)} \Norm{\w}_{L^{1}(\Omega)} \,,
\]
and analogous estimates hold for $\wrho^n(x):= \int_\Omega w(|x-z|) \varrho(z) \d z$.
Consequently, we have
\[
    |u^n_\eps(x) - u^n_\eps(y)| &\leq&
              \frac{1}{\eps} |\wj^n(x)-\wj^n(y)|
            + \frac{1}{\eps^2} |\wj^n(x)-\wj^n(y)| \wrho^n(y)
            + \frac{1}{\eps^2} |\wrho^n(x)-\wrho^n(y)| |\wj^n(y)| \\
    &\leq& \frac{1}{\eps} \Norm{j^n}_{L^\infty(\Omega)} \left(
              1 + \frac{2}{\eps} \Norm{\varrho^n}_{L^\infty(\Omega)} \Norm{\w}_{L^{1}(\Omega)} \right) 
         \omega(x-y) \,.
\]
This uniform equicontinuity together with the uniform boundedness $|u^n_\eps|\leq 1$ allows us
to apply the Arzela-Ascoli Lemma and obtain the compactness of the Schauder operator.
Therefore, for every fixed $\eps>0$, we have a nonnegative solution
$(p^+_\eps,p^-_\eps)$ of the regularized system~\fref{kinetic1}, \fref{kinetic2}, \fref{u_eps},
uniformly bounded (with respect to $\eps$) in $L^\infty([0,T]\times \Omega)$.

\medskip\noindent\textbf{Step 3.}
Finally, we pass to the limit $\eps\to 0$.
Due to the uniform boundedness,
a subsequence of $p^\pm_\eps\rightharpoonup p^\pm$ weakly* in
$L^\infty([0,T]; L^\infty(\Omega))$;
we need to show that the limit of the nonlinear terms $\varrho_\eps u_\eps$ and $u_\eps^2 j_\eps$
is $\varrho u$ and, resp., $u^2 j$, with $u$ given by~\fref{u_smooth}, or, equivalently, \fref{hydrodynamic3}.
The limit passage in the distributional formulation of the term $\varrho_\eps u_\eps$
with a test function $\phi\in C^\infty_c([0,T)\times\Omega)$ is performed as follows:
\[
    \left| \int_0^\infty \int_\Omega (\varrho_\eps u_\eps - \varrho u) \phi \d x\d t \right| \leq
       \left| \int_0^\infty \int_\Omega (\varrho_\eps -\varrho) u\phi \d x\d t \right| +
       \left| \int_0^\infty \int_\Omega \varrho_\eps (u_\eps- u) \phi \d x \d t\right| \,.
\]
The first term vanishes in the limit $\eps\to 0$ due to the weak* convergence
of $\varrho_\eps$ towards $\varrho$ with $u\phi$ a valid test function.
Concerning the second term, we will show that $\int_\Omega \varrho_\eps (u_\eps- u) \phi \d x$
tends to zero for almost all $t\in [0,T]$ and conclude the convergence of the time integral
by the Lebesgue dominated convergence theorem.
Let us fix $t\in [0,T]$ and define the set $\mathcal{S}_\delta$ by
\[
    \mathcal{S}_\delta := \left\{x\in\Omega;\, \wrho(x) < \delta \right\} \,,
\]
with $\delta\geq 0$ and $\wrho$ defined in~\fref{u_eps}.
Then, we have
\[
    \left| \int_{\mathcal{S}_\delta} \varrho_\eps (u_\eps- u) \phi \d x\right|
     &\leq&   2 \int_{\mathcal{S}_\delta} \varrho_\eps |\phi| \d x
     \;\mathop{-\!\!\!\longrightarrow}_{\!\eps\to 0}\;   2 \int_{\mathcal{S}_\delta} \varrho |\phi|  \d x \\
     &=&    2  \int_{\mathcal{S}_\delta\setminus \mathcal{S}_0} \varrho |\phi| \d x
     \;\mathop{-\!\!\!\longrightarrow}_{\!\delta\to 0}\; 0  \,,
\]
where the second line is due to $\varrho=0$ on $\mathcal{S}_0$
and because $\mbox{meas}(\mathcal{S}_\delta\setminus \mathcal{S}_0)$ tends to zero as $\delta\to 0$.
Next, for $x\in \Omega\setminus\mathcal{S}_\delta$,
\[
    |u_\eps(x) - u(x)| &\leq& \frac{|\wj_\eps(x)| |\wrho_\eps(x)-\wrho(x)|}{\wrho(x)\wrho_\eps(x)}
         + \frac{|\wj_\eps(x)-\wj(x)|}{\wrho(x)}  \\
        &\leq&  \frac{1}{\delta} \left( |\wrho_\eps(x)-\wrho(x)| + |\wj_\eps(x)-\wj(x)| \right) \,,
\]
where $\wj_\eps(x) = \int_\Omega w(|x-z|) j_\eps(z) \d z$ and $\wrho_\eps(x) = \int_\Omega w(|x-z|) \varrho_\eps(z) \d z$.
Therefore, due to the uniform convergence of $\wrho_\eps$ and $\wj_\eps$ to $\wrho$ and, resp., $\wj$ on $\Omega$
(implied by the uniform equicontinuity and boundedness of the families $\{\wrho_\eps\}_{\eps>0}$
and $\{\wj_\eps\}_{\eps>0}$),
we have
\[
    \left| \int_{\Omega\setminus \mathcal{S}_\delta} \varrho_\eps (u_\eps- u) \phi \d x\right|
    &\leq& \left( \mbox{sup}_{\Omega\setminus \mathcal{S}_\delta} |u_\eps-u| \right)
           \int_\Omega \varrho_\eps |\phi| \d x \\
    &\leq& \frac{1}{\delta} \left( \Norm{\wrho_\eps-\wrho}_{L^\infty(\Omega)} +
             \Norm{\wj_\eps-\wj}_{L^\infty(\Omega)} \right)
        \Norm{\varrho_\eps}_{L^\infty(\Omega)} \Norm{\phi}_{L^1(\Omega)} 
    \;\mathop{-\!\!\!\longrightarrow}_{\!\eps\to 0}\; 0 \,.
\]
We conclude by passing $\delta\to 0$.
The limit passage in the term $u_\eps^2 j_\eps$ is performed similarly
(note that $|j_\eps| \leq \varrho_\eps$).
\endproof

\noindent
It is worth noting that the assumption~\ref{A1} of Theorem~\ref{thm:existence} can be relaxed.
In fact, we posed the requirement [A1] of boundedness of $\w$ on $\Omega$
in order to establish the definition~\fref{u_i^loc} of $u_i^{loc}$
in the discrete model and pass to the limit $N\to\infty$.
However, at the level of the kinetic or hydrodynamic description,
we may relax this to $\w\in L^1(\Omega)$.
Moreover, formally it is possible to consider even singular weights,
in particular $\w = \delta_0$, which leads to $u = j/\varrho$ and removes the nonlocality.
In fact, one can see the choice $\w=\delta_0$ as the limiting case
when the interaction radius shrinks to zero:
for almost all $x\in\Omega$ such that $\varrho(x)\neq 0$, one has
\[
   \lim_{\sigma\to 0}
   \frac{\int_\Omega \chi_{[0,\sigma]}(|x-z|) j(z) \d z}{\int_\Omega \chi_{[0,\sigma]}(|x-z|) \varrho(z) \d z}
           = \frac{j(x)}{\varrho(x)} \,,
\]
where $\chi_{[0,\sigma]}$ is the characteristic function of the interval $[0,\sigma]$.
One can interpret this as a model where only pointwise local observations of the system
are possible.

By a slight modification of the proof of Theorem~\ref{thm:existence}
it is possible to show that given a sequence of weights $\w_n$ converging strongly in $L^1(\Omega)$
to $\w$,
the solutions $(\rho_n, j_n)$ corresponding to the weights $\w_n$ converge weakly*
in $L^\infty([0,T]; L^\infty(\Omega))$ to the solution corresponding to the weight $\w$.
However, we need the condition [A2'] to be satisfied uniformly;
consequently, the question whether and how
the solution corresponding to $\w=\delta_0$ can be obtained
as a limit of solutions with $\w_\sigma = \chi_{[0,\sigma]}$ as $\sigma\to 0$
remains an interesting open problem.
An even more interesting question is what is the limit of the discrete model
as $N\to\infty$ if the interaction radius is shrinking as some power of $1/N$.
This question is studied in Section~\ref{sec:Shrinking} below.

\section{Long time behaviour}\label{sec:Long time}
\begin{figure}\label{fig:num1}
{\centering \begin{tabular}[h]{ccc}
\resizebox*{0.33\linewidth}{!}{\includegraphics{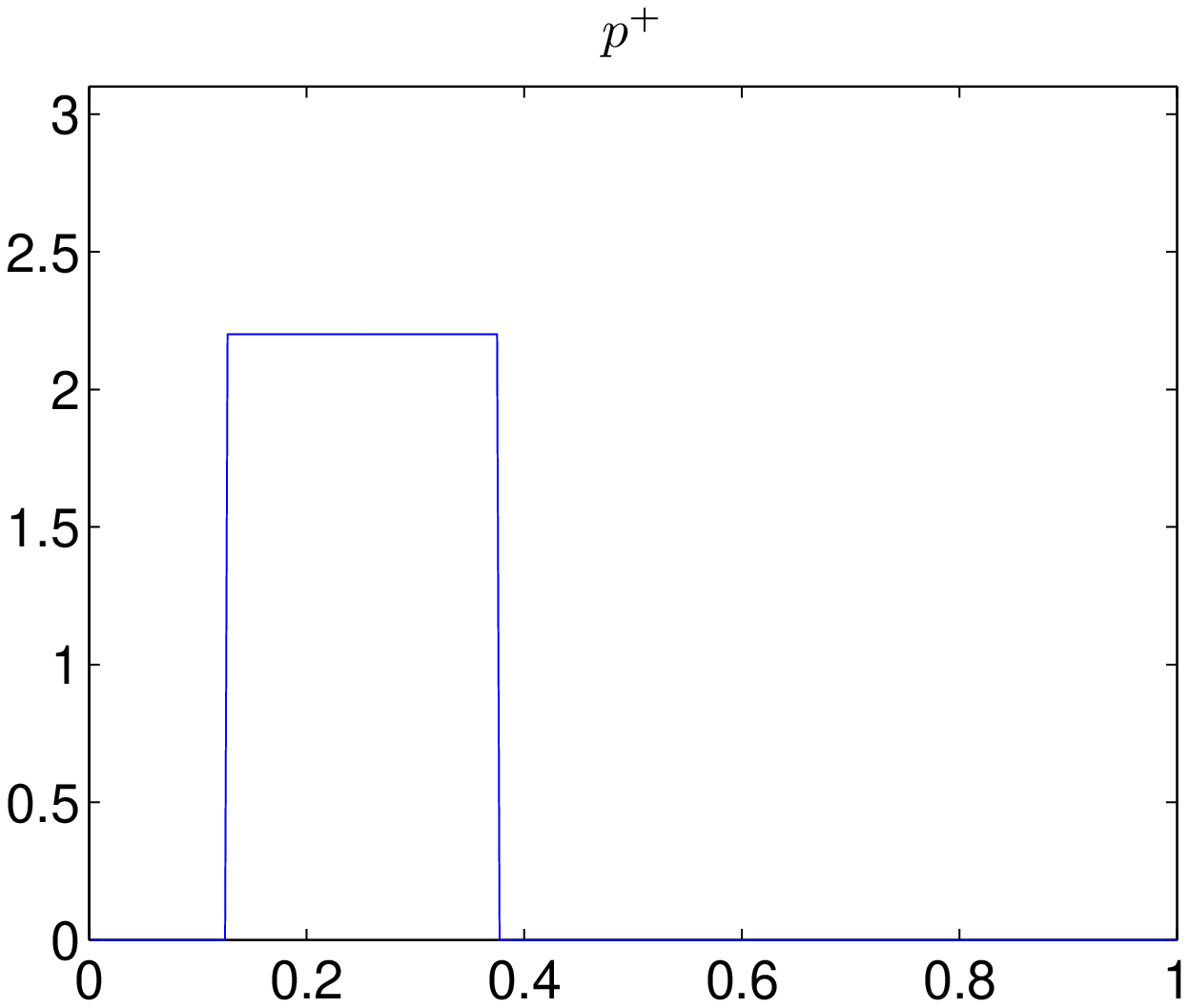}} &
\resizebox*{0.33\linewidth}{!}{\includegraphics{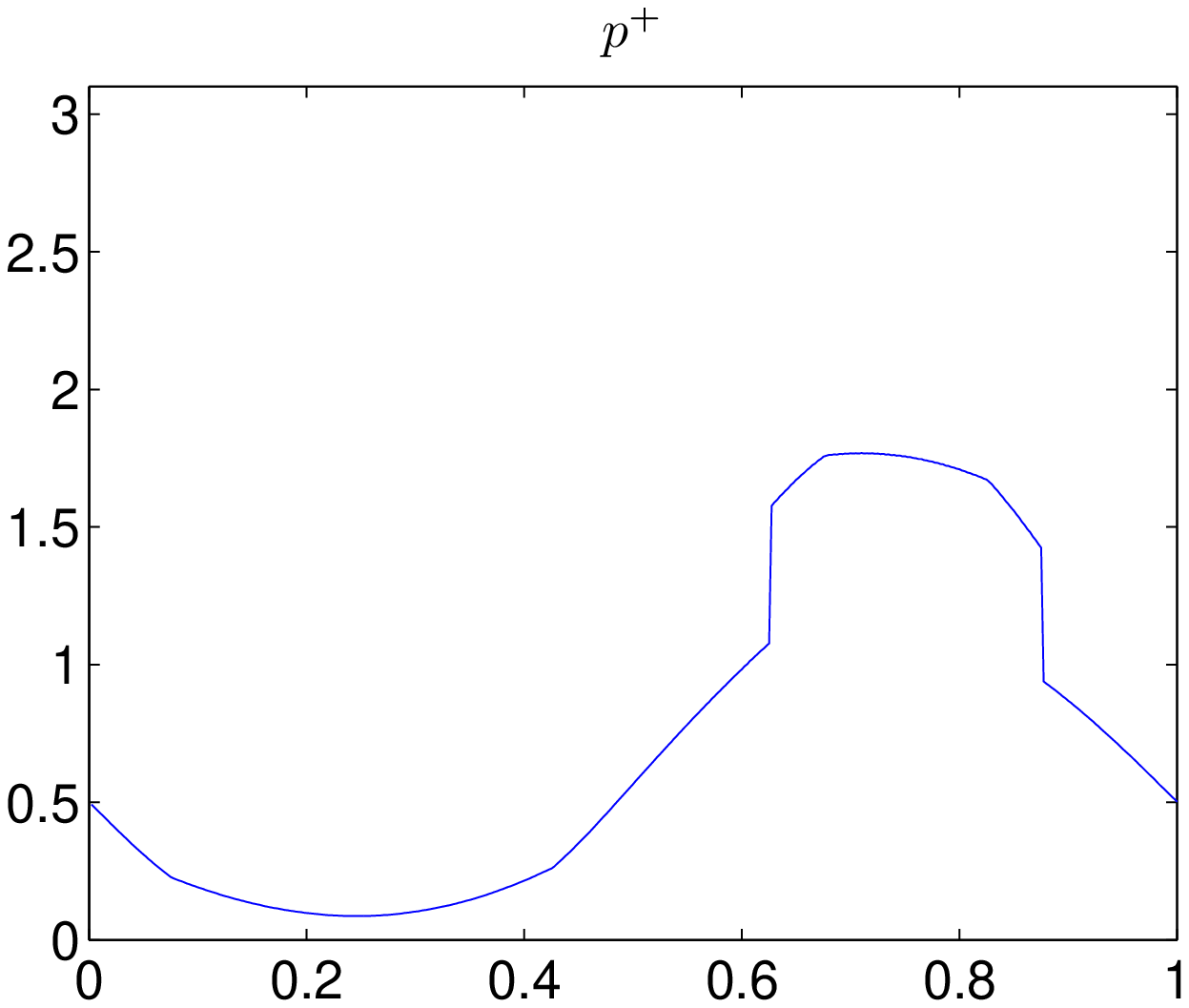}} &
\resizebox*{0.33\linewidth}{!}{\includegraphics{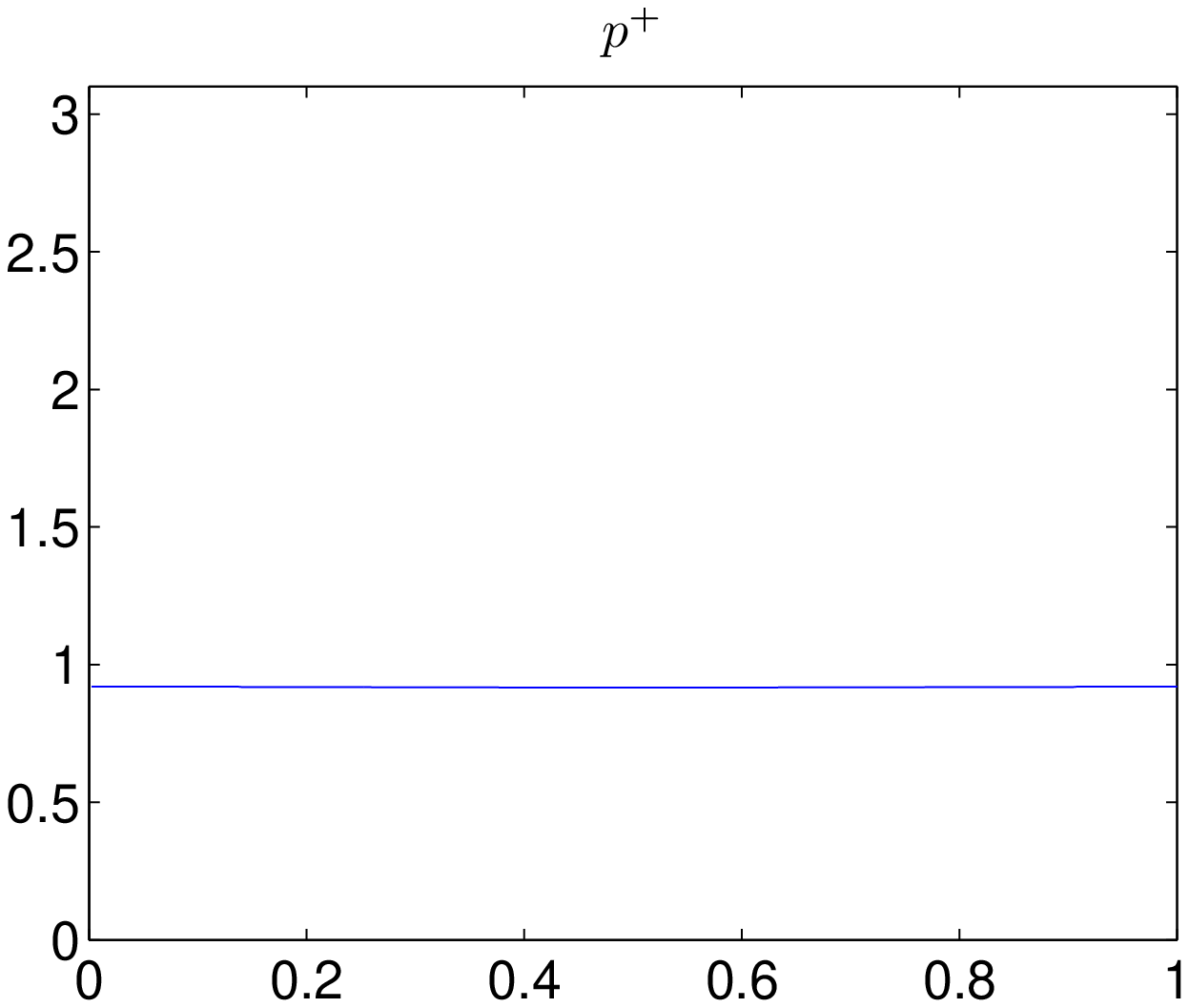}} \\
\resizebox*{0.33\linewidth}{!}{\includegraphics{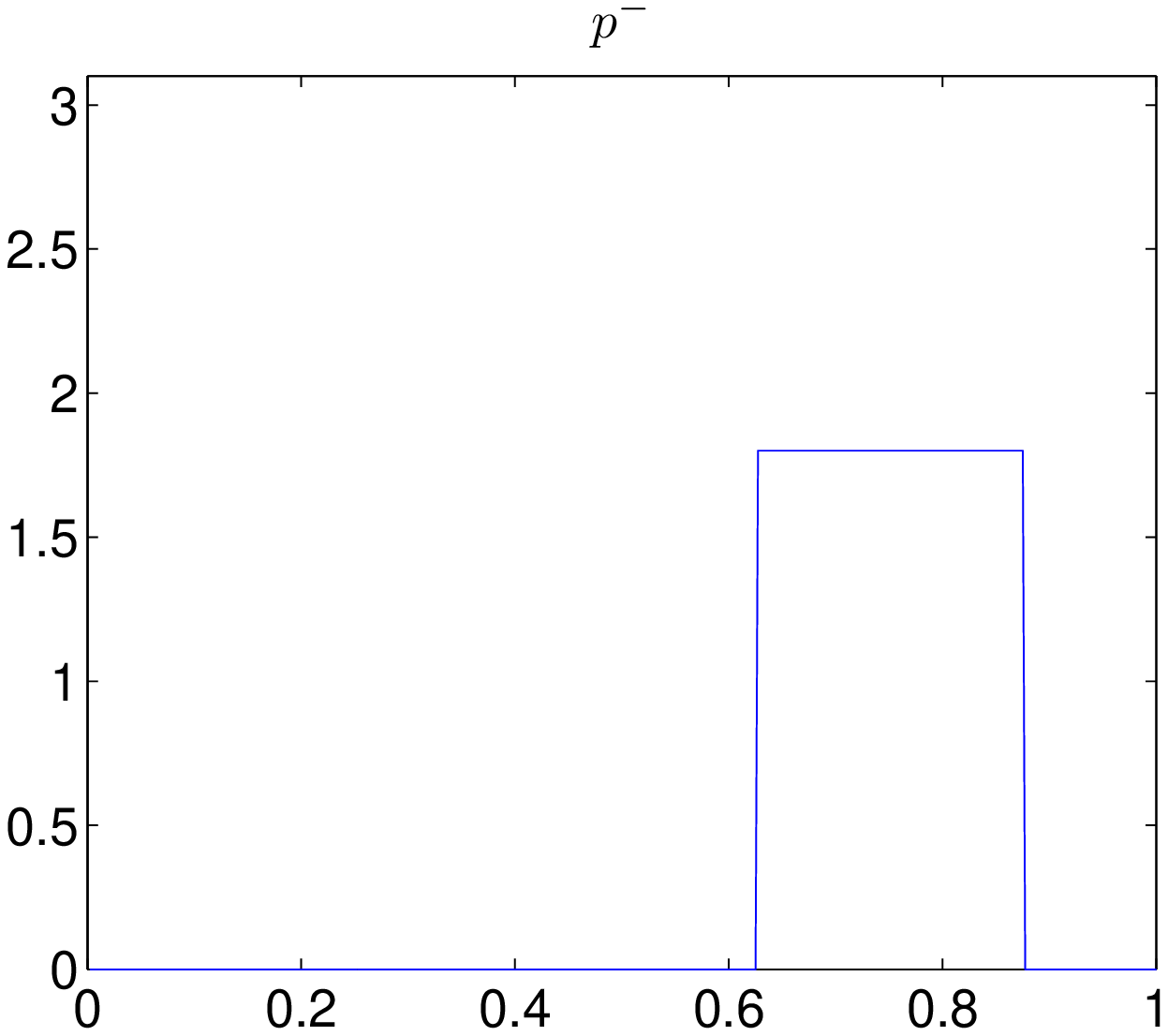}} &
\resizebox*{0.33\linewidth}{!}{\includegraphics{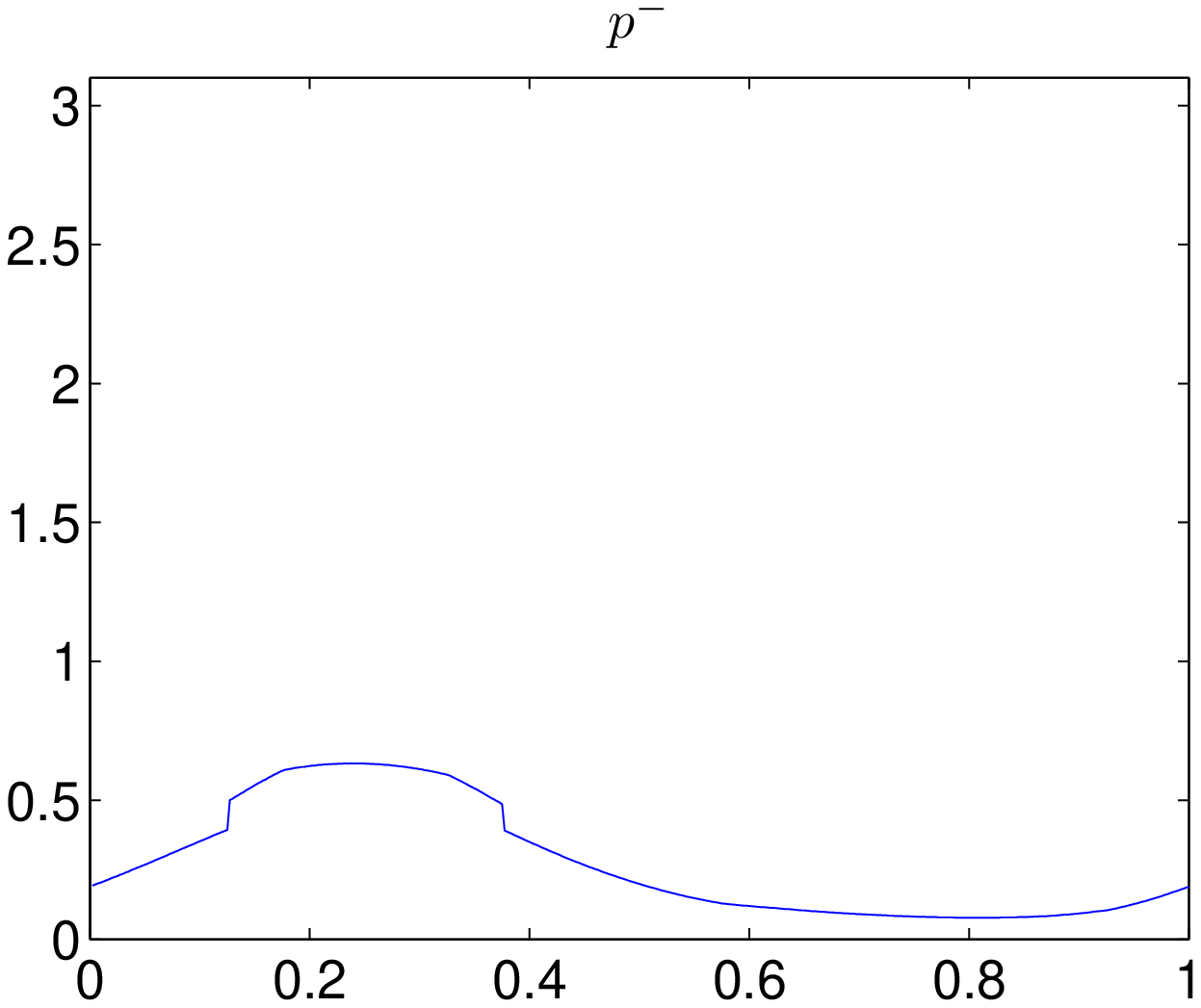}} &
\resizebox*{0.33\linewidth}{!}{\includegraphics{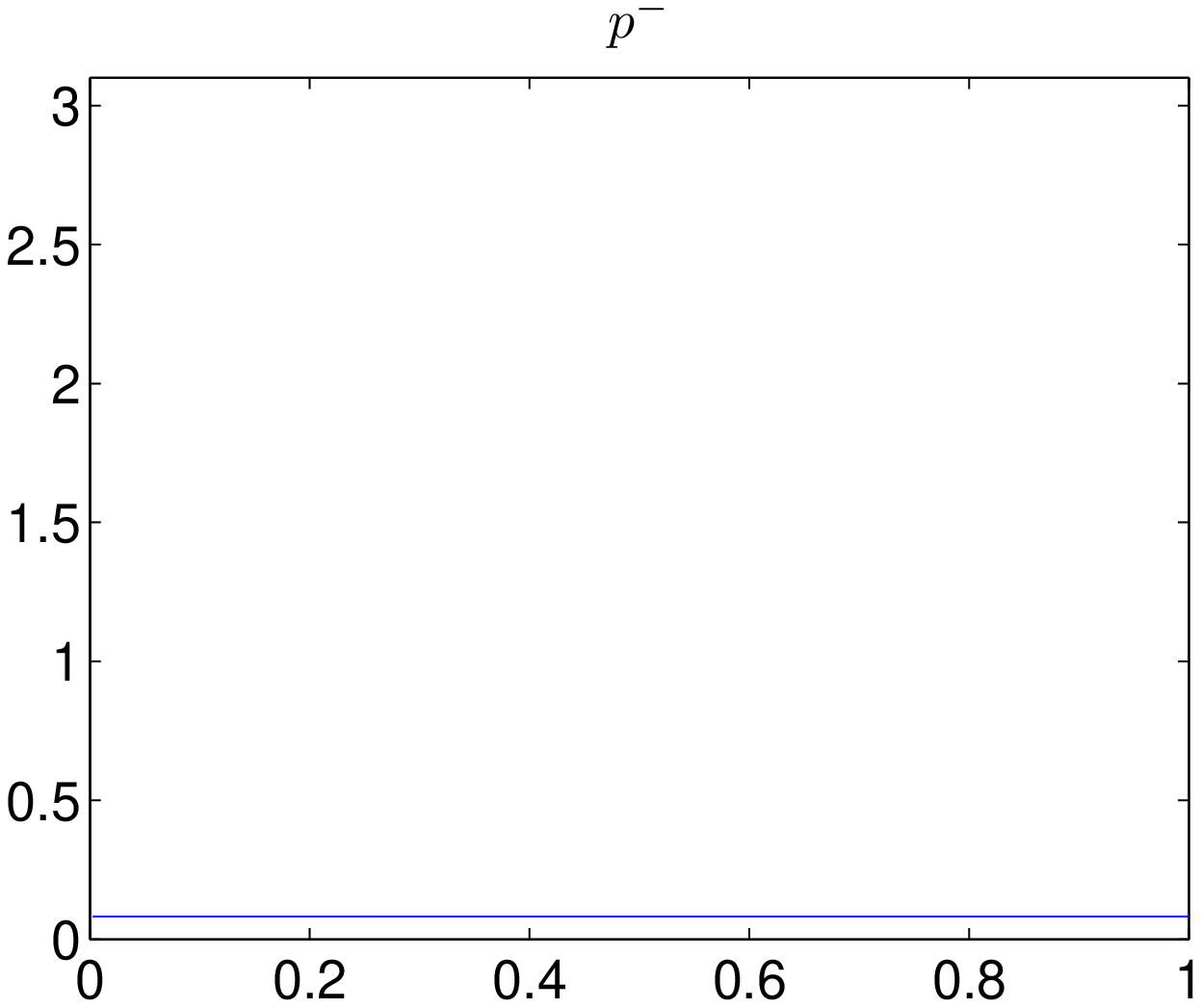}} \\
\resizebox*{0.33\linewidth}{!}{\includegraphics{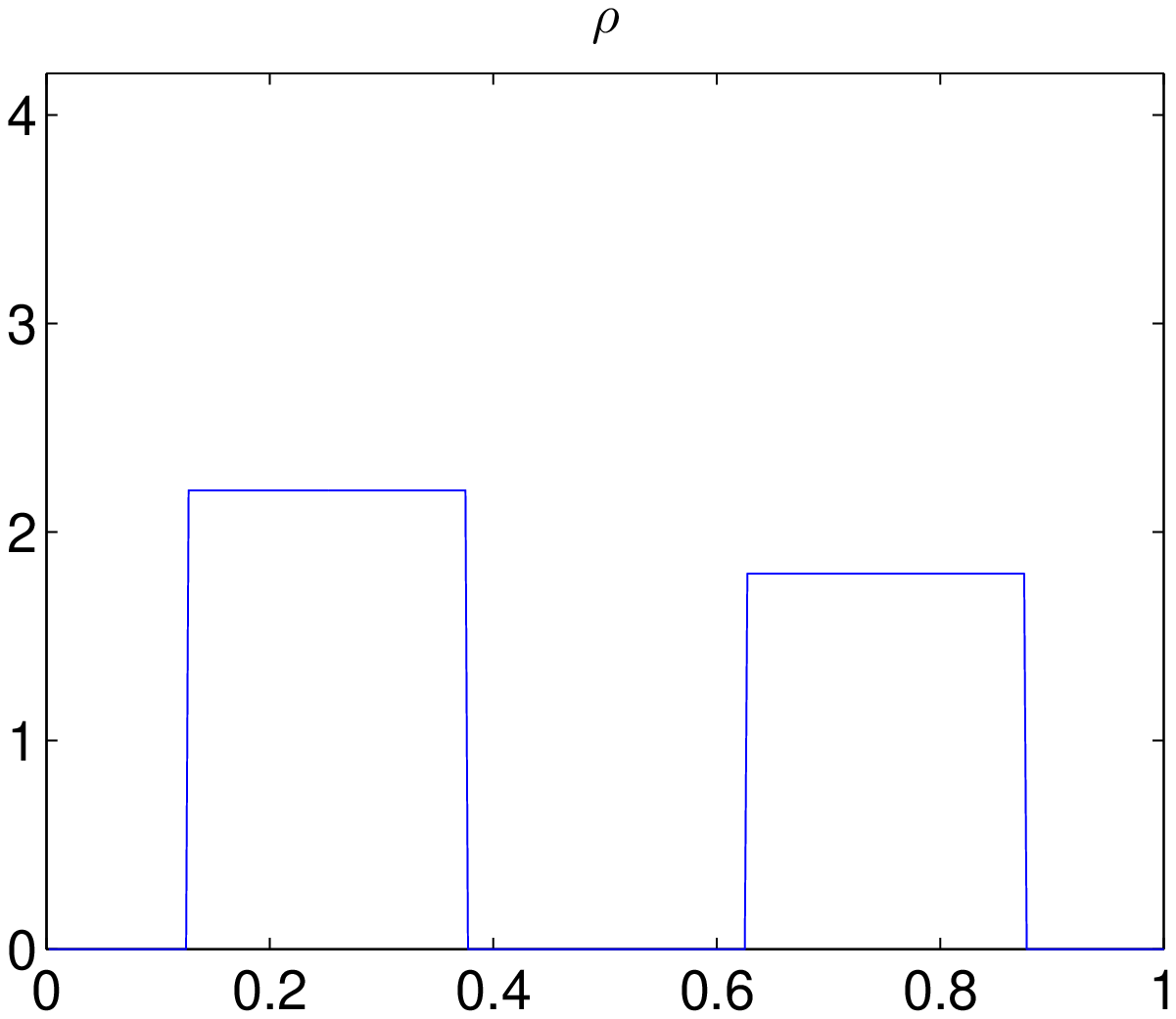}} &
\resizebox*{0.33\linewidth}{!}{\includegraphics{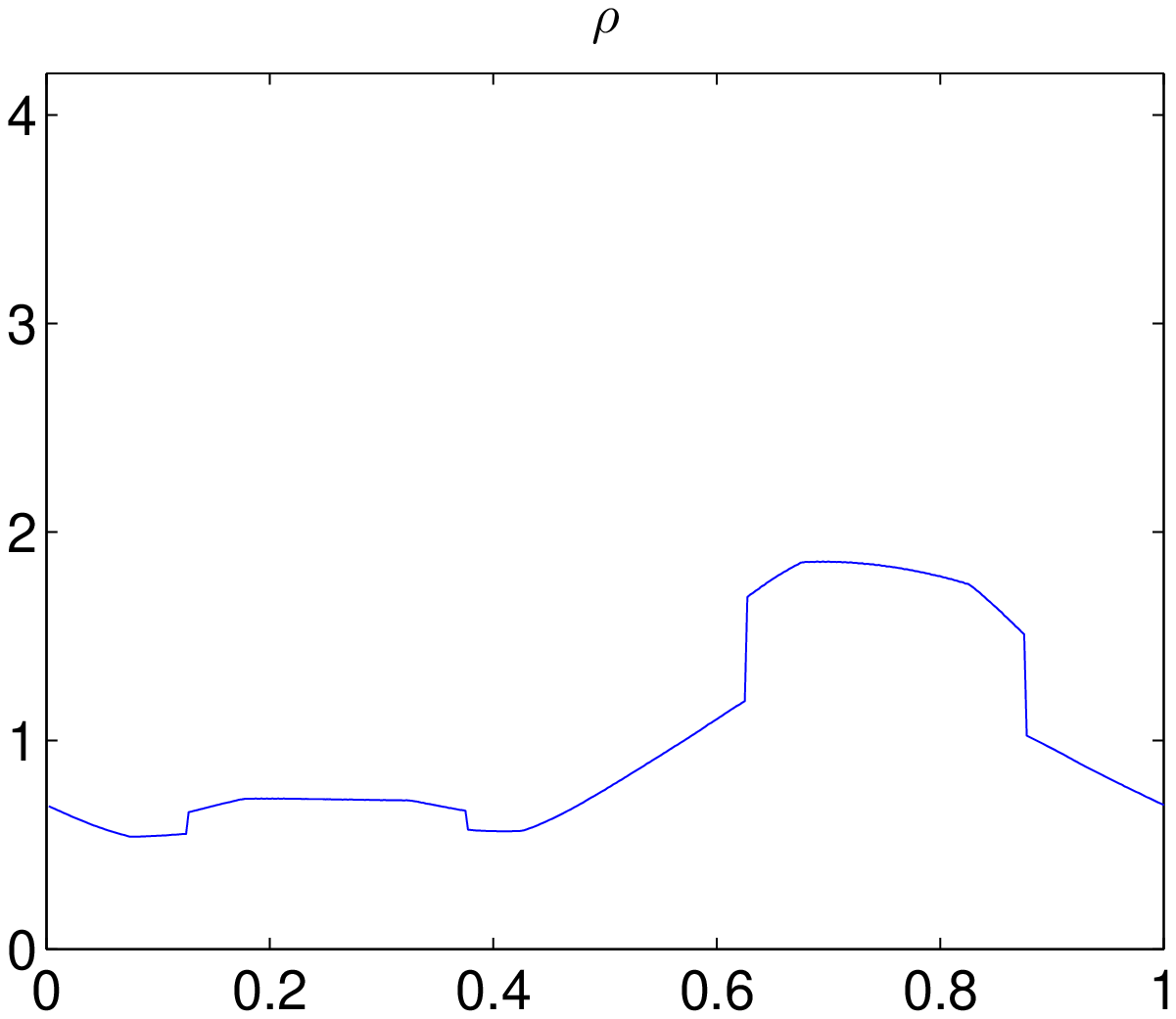}} &
\resizebox*{0.33\linewidth}{!}{\includegraphics{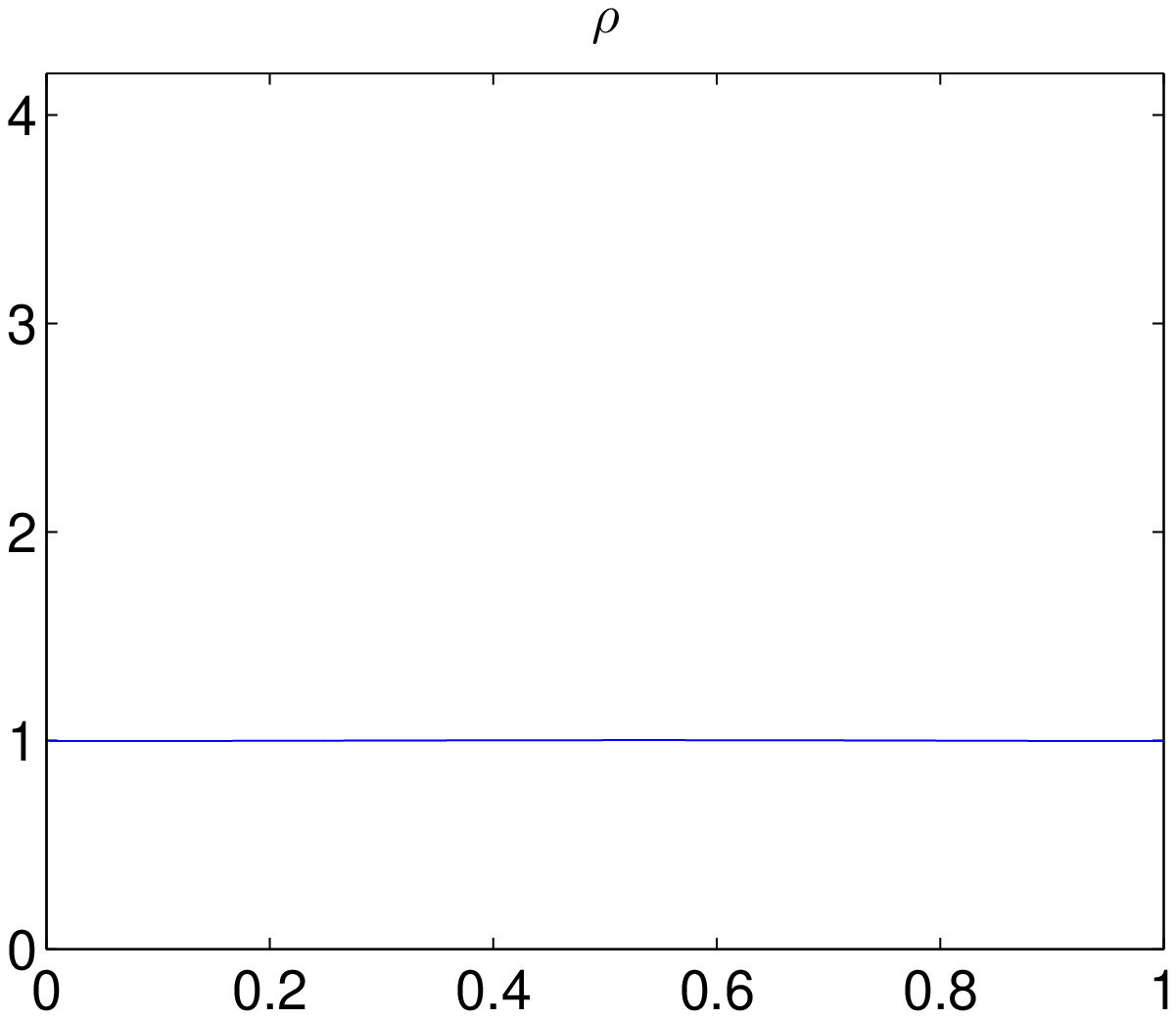}} \\
\resizebox*{0.33\linewidth}{!}{\includegraphics{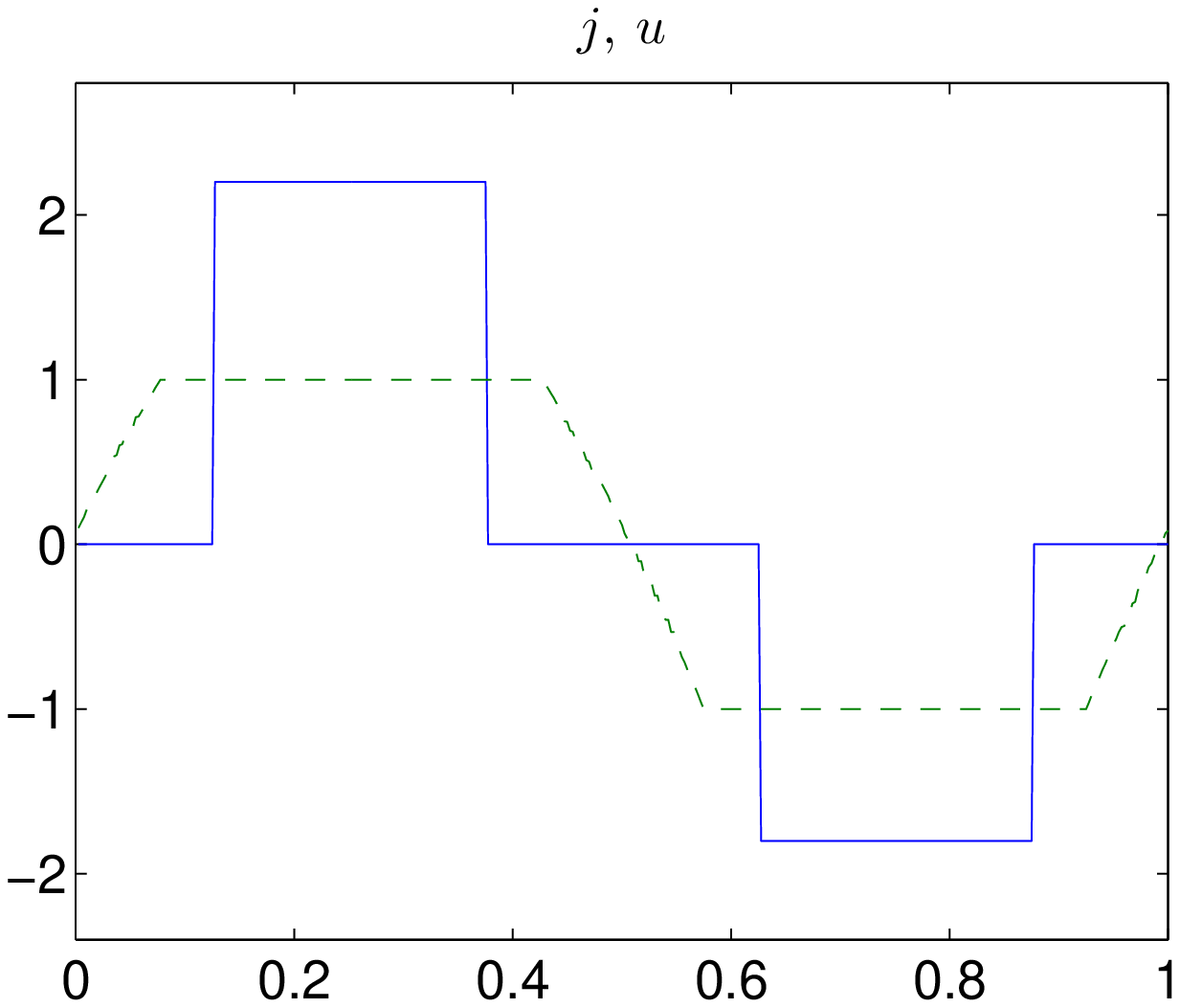}} &
\resizebox*{0.33\linewidth}{!}{\includegraphics{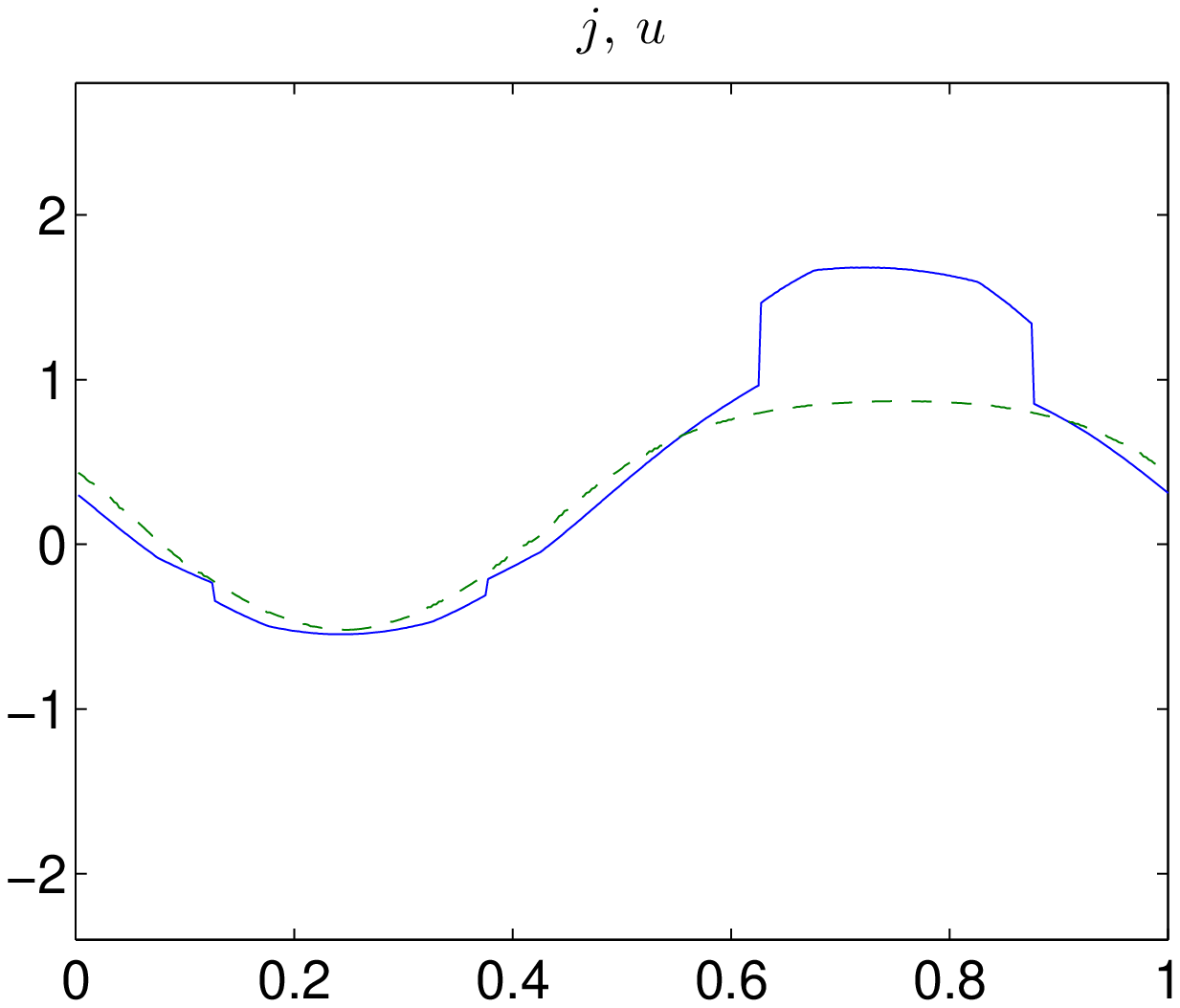}} &
\resizebox*{0.33\linewidth}{!}{\includegraphics{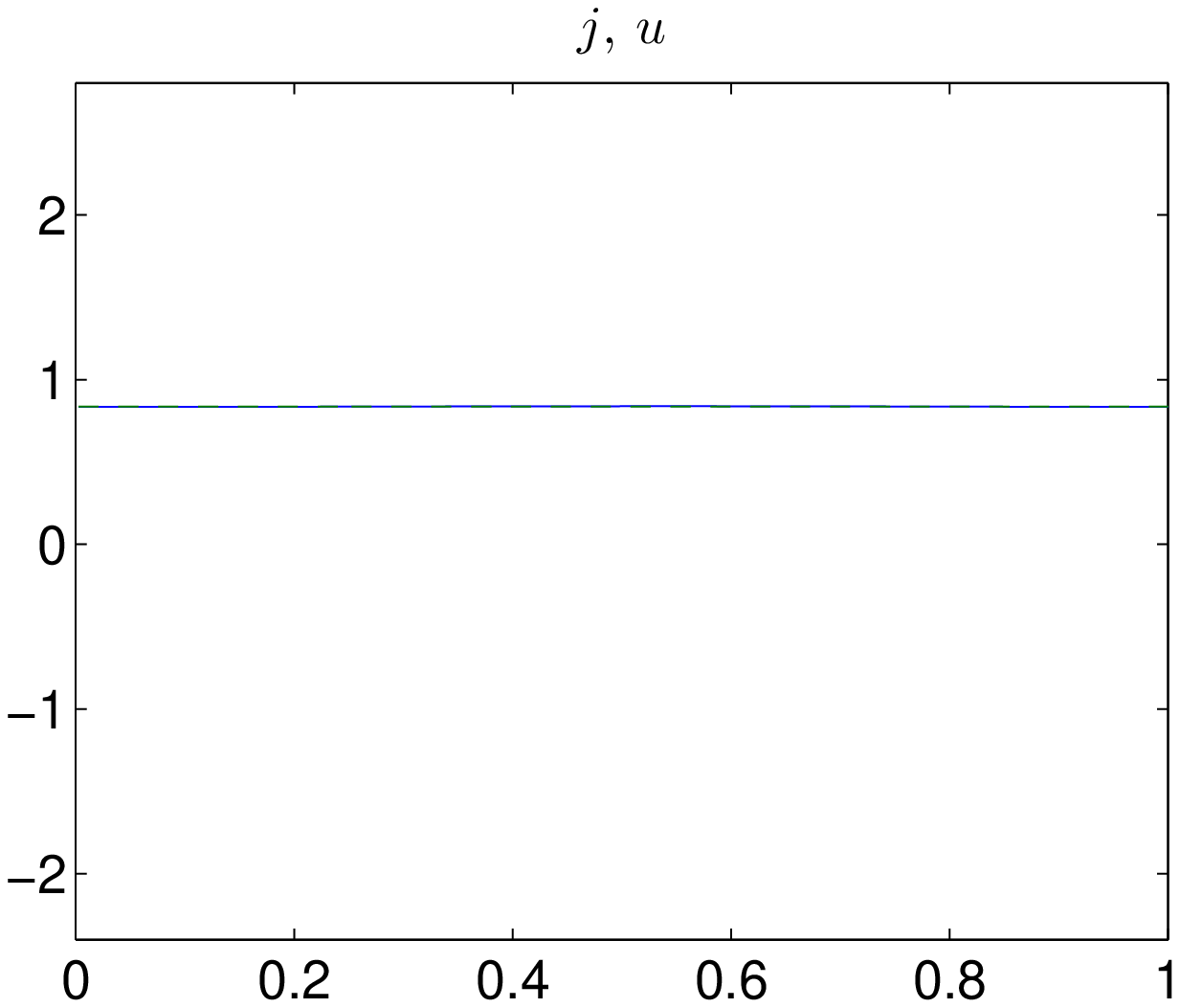}} \\
$t=0.0$ & $t=2.5$ & $t=25.0$
\end{tabular}\par}
\caption{Numerical results with $\w=\chi_{[0,0.2]}$, $b=1$ and $\gamma_0=0.3$:
$p^+$ and $p^-$ converge to constant states, the mass density $\varrho$ converges to $1$,
the flux $j$ (full line) and averaged velocity $u$ (dashed line) converge to a constant.
}
\end{figure}

In this Section we provide several conjectures about the long time behaviour of the kinetic system~\fref{u_smooth}, \fref{kinetic1}, \fref{kinetic2} or,
equivalently, the hydrodynamic system~\fref{hydrodynamic1}--\fref{hydrodynamic3}.
To get some insight into the long time dynamics, we start with a numerical example.
We solve the kinetic system using standard semi-implicit finite difference method with upwinding.
The initial condition is $p^+_0 = 2.2$ on $[0.125, 0.375]$ and zero otherwise, $p^-_0=1.8$ on $[0.625, 0.875]$
and zero otherwise. In Figure~\ref{fig:num1} we show the results for the choice of parameters
$b=1$ and $\gamma_0=0.3$ and the weight function $\w=\chi_{[0,0.2]}$.

\noindent
\vspace{0.3cm}

We conjecture that with regular weights $\w$ satisfying [A1] and [A2'], the solutions $(\rho, j)$ 
to~\fref{hydrodynamic1}--\fref{hydrodynamic3} converge to the constants $\rho\equiv 1$
and $j\equiv j_s$ with some $j_s\in\R$, exponentially fast as $t\to\infty$.
Moreover, in the large noise case $\gamma_0 > b$, we hypothesize that $j_s \equiv 0$.
Unfortunately, we are only able to provide an analytical proof
in the rather special case $\w\equiv 1$ and $\gamma_0 > b$:

\begin{lemma}\label{lemma:long_time_beh_reg}
Assuming $\w\equiv 1$, we have $u(t,x) \equiv u(t)$, where $u(t)$ satisfies
the ordinary differential equation
\(    \label{ODE_u}
    \dot u = - 2(\gamma_0 + b(u^2-1)) u \,,
\)
subject to the initial condition $u(0) = \int_\Omega j(0,x) \d x$.
Moreover,
\begin{enumerate}[leftmargin = 1.0 cm, label = {(}\roman*{)}]
\item if $\gamma_0 \leq b$, then $\lim_{t\to\infty} u(t) = \sign(u(0)) \sqrt{1-\gamma_0 b^{-1}}$,
\item if $\gamma_0 > b$, then $\lim_{t\to\infty} u(t) = 0$ and
\(   \label{est_u}
    |u(t)| \leq |u(0)| e^{-2(\gamma_0 - b)t} \,.
\)
Moreover, $j$ converges to zero exponentially fast in the $L^2$-sense:
\[
    \int_\Omega j^2(t,x) \d x \leq c e^{-4\gamma_0 t}
\]
for a suitable constant $c$.
\end{enumerate}
\end{lemma}

\begproof
Integrating~\fref{hydrodynamic2} with $\w\equiv 1$ with respect to $x\in\Omega$, we obtain~\fref{ODE_u}.
This has the stationary state $u=0$, which is stable if and only if $\gamma_0 > b$.
Moreover, if $\gamma_0 < b$, two additional stable stationary states
$u = \pm\sqrt{1-\gamma_0 b^{-1}}$ exist.
This establishes the first statement.
Further, we have
\[
    \tot{}{t} \left(\frac12 u^2\right) = -2 (\gamma_0 + b(u^2-1)) u^2 \leq -2 (\gamma_0 - b) u^2 \,,
\]
and an application of the Gronwall lemma gives~\fref{est_u}.

To prove the convergence of $j$ to zero, we consider the identity
\(    \label{estimate}
    \frac12 \tot{}{t} \int_\Omega \left( \varrho^2 + j^2 \right) \d x = -2 (\gamma_0 + b(1+u^2)) \int_\Omega j^2 \d x
        + 4b u \int_\Omega \varrho j \d x \,.
\)
An application of the Cauchy-Schwartz inequality and the decay rate~\fref{est_u} yield
\[
    4b u \int_\Omega \varrho j \d x \leq 2b \int_\Omega j^2 \d x + 2b u^2 \int_\Omega \varrho^2 \d x
      \leq 2b \int_\Omega j^2 \d x + 2b |u(0)|^2 e^{-4(\gamma_0 - b)t}\int_\Omega \varrho^2 \d x \,.
\]
Then, an integration of~\fref{estimate} in time leads to
\[
    \int_\Omega \varrho^2 (T,x) \d x \leq c_0 + 4b |u(0)|^2 e^{-4(\gamma_0 - b)t}\int_0^T \int_\Omega \varrho^2(t,x) \d x\d t \,,
\]
with $c_0 := \int_\Omega \varrho_0^2(x) + j_0^2(x) \d x$.
Consequently, by the Gronwall lemma, $\int_\Omega \varrho^2 \d x$ is bounded uniformly in time by a constant
$c_1$ if $\gamma_0 > b$.
Inserting this information into~\fref{estimate} gives
\[
    \int_\Omega j^2 (T,x) \d x &\leq& c_0 - 4\gamma_0 \int_0^T \int_\Omega j^2(t,x) \d x\d t 
       + 4b c_1 |u(0)|^2 \int_0^T e^{-4(\gamma_0 - b) t} \d t \\
   &\leq& c - 4\gamma_0 \int_0^T \int_\Omega j^2(t,x) \d x\d t 
\]
and an application of the Gronwall lemma yields the second statement.
\endproof


\subsection{The case $\w=\delta_0$}\label{subsect:long_time_beh_sing}
\begin{figure}\label{fig:num2}
{\centering \begin{tabular}[h]{ccc}
\resizebox*{0.33\linewidth}{!}{\includegraphics{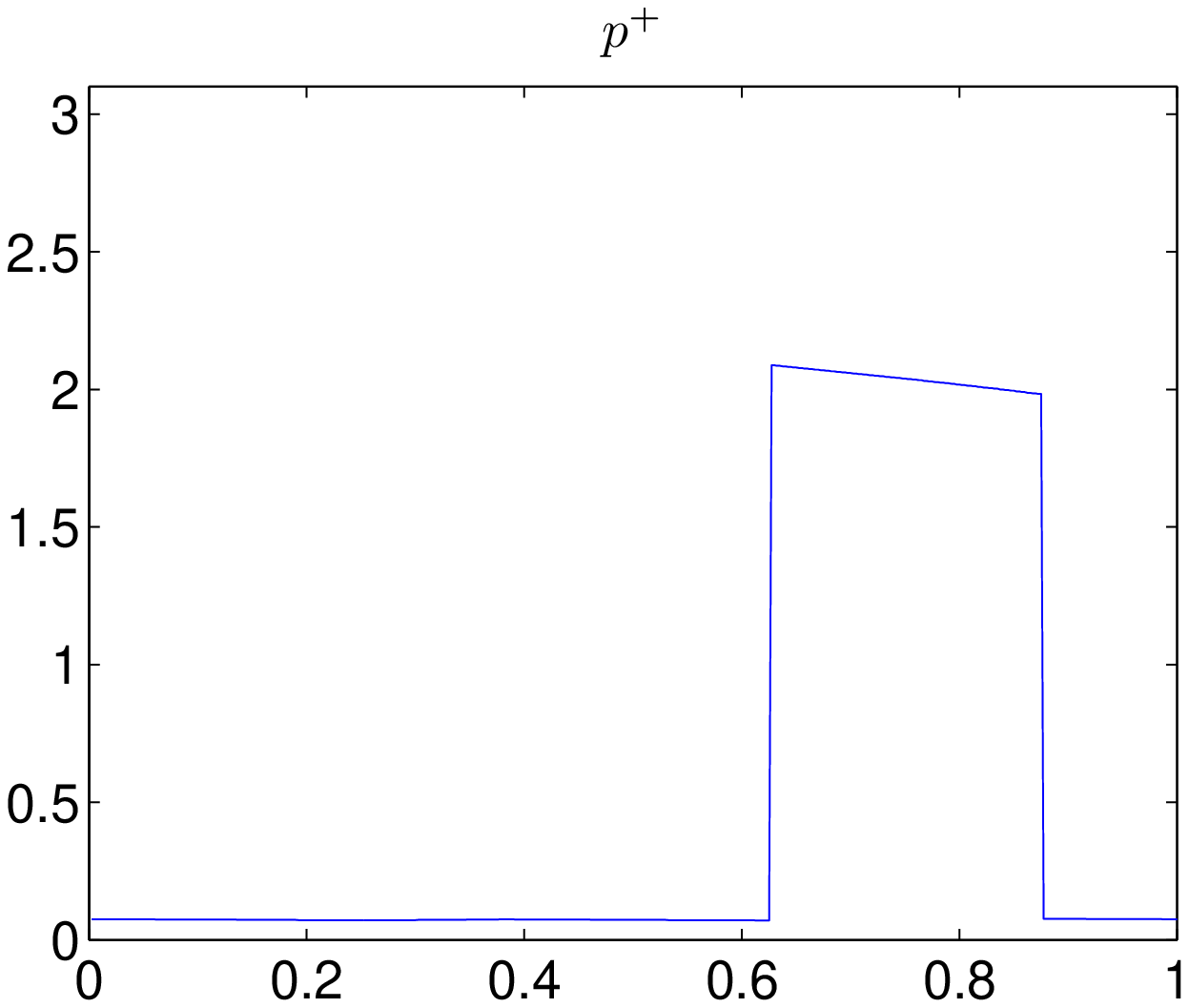}} &
\resizebox*{0.33\linewidth}{!}{\includegraphics{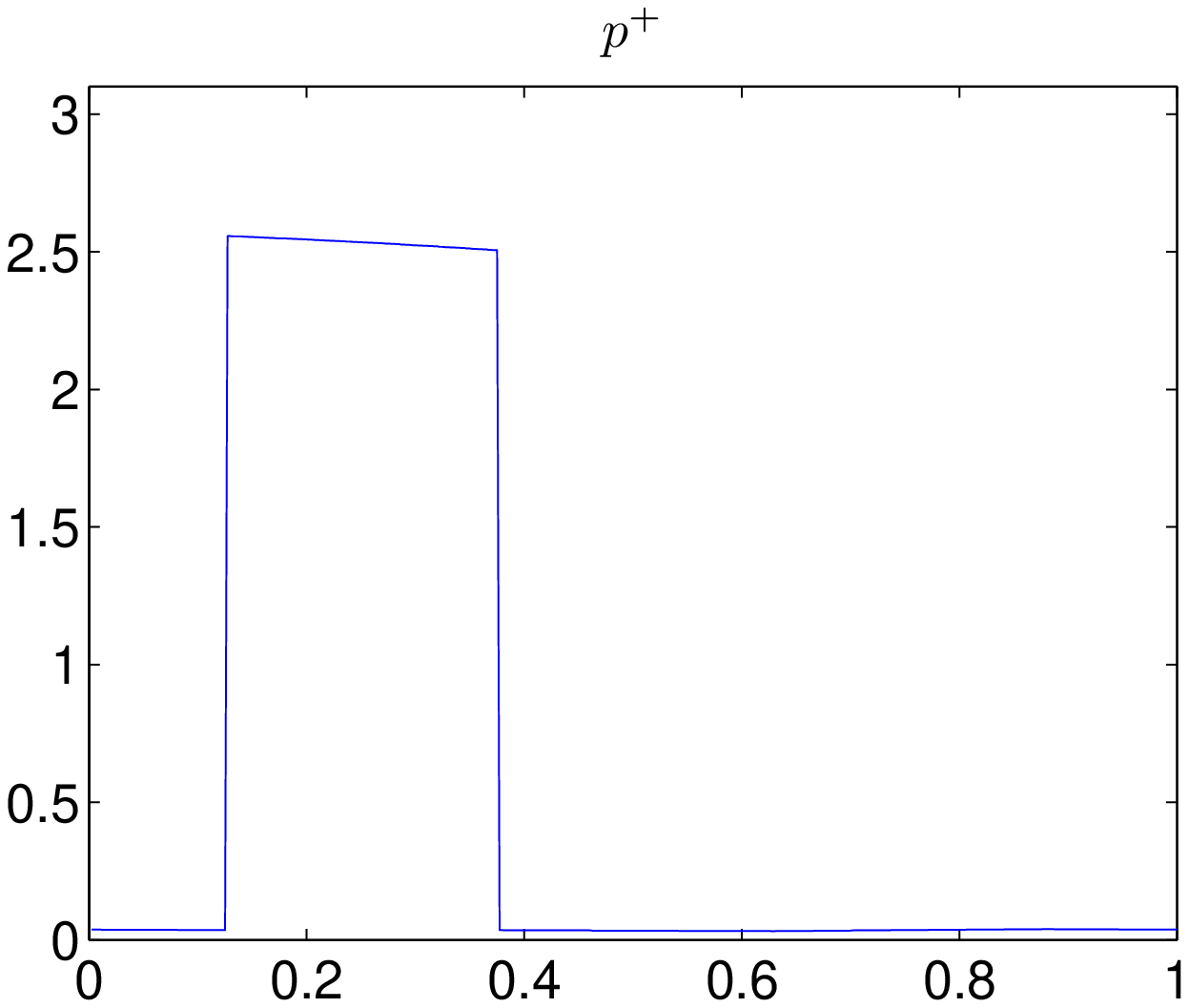}} &
\resizebox*{0.33\linewidth}{!}{\includegraphics{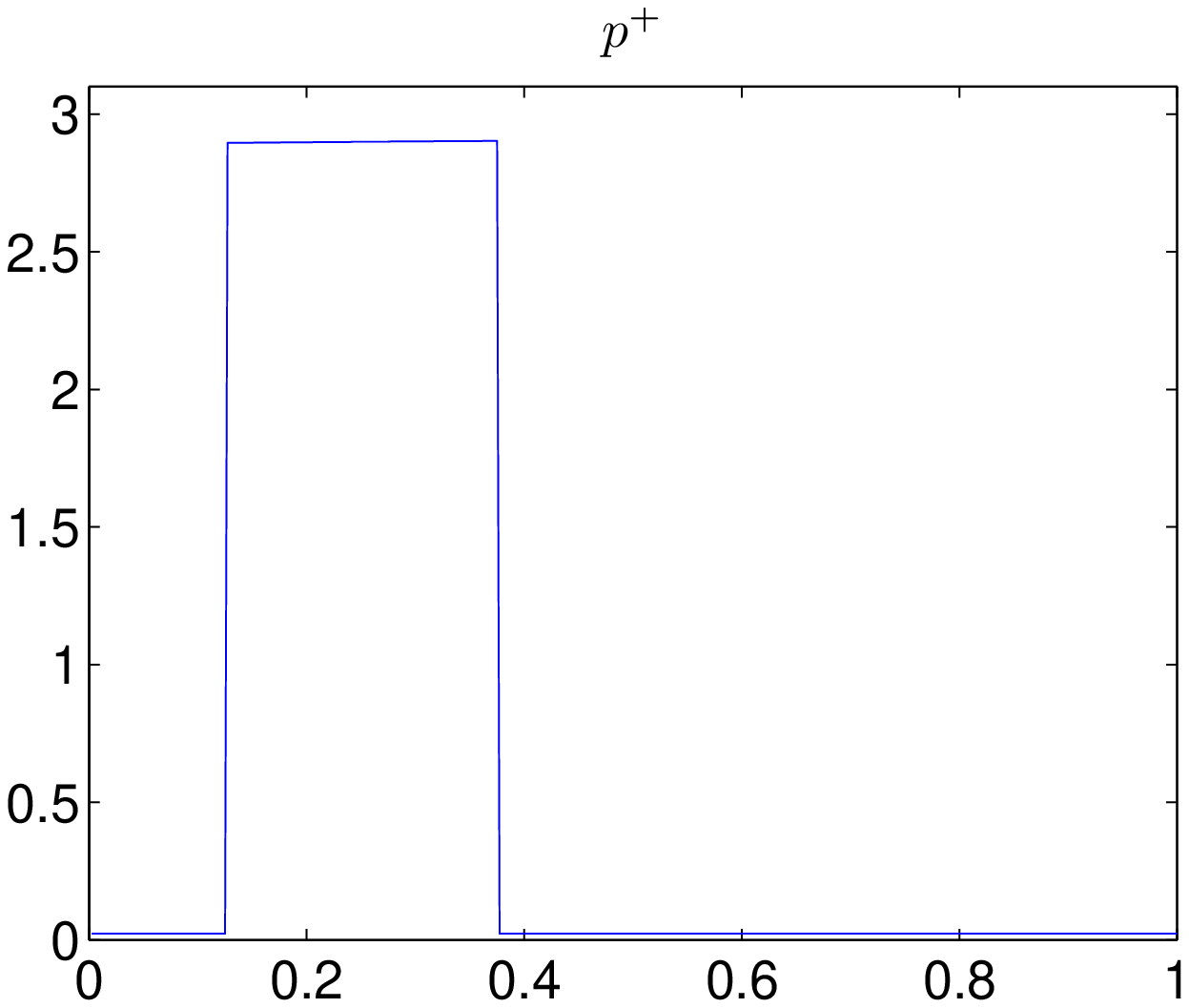}} \\
\resizebox*{0.33\linewidth}{!}{\includegraphics{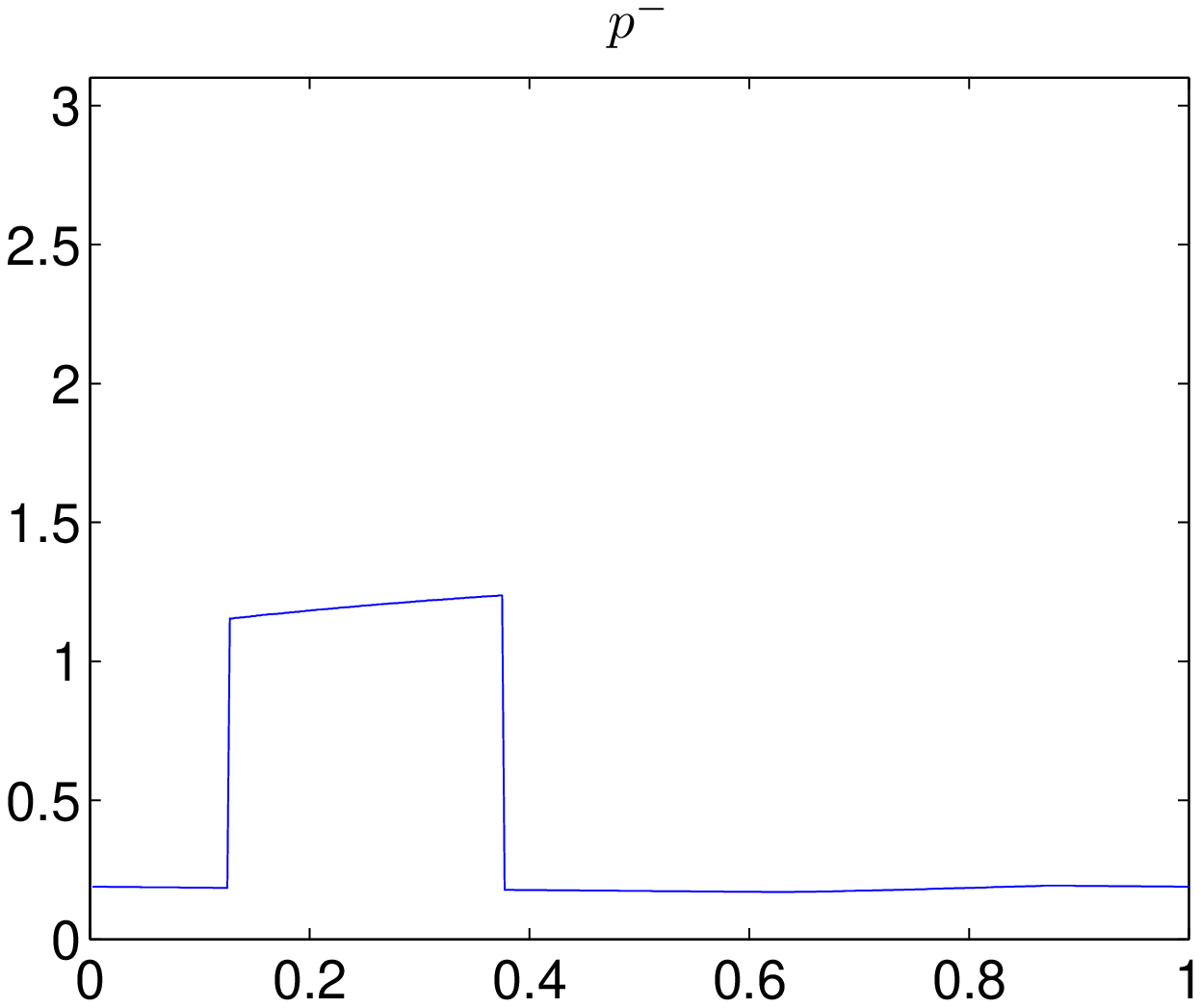}} &
\resizebox*{0.33\linewidth}{!}{\includegraphics{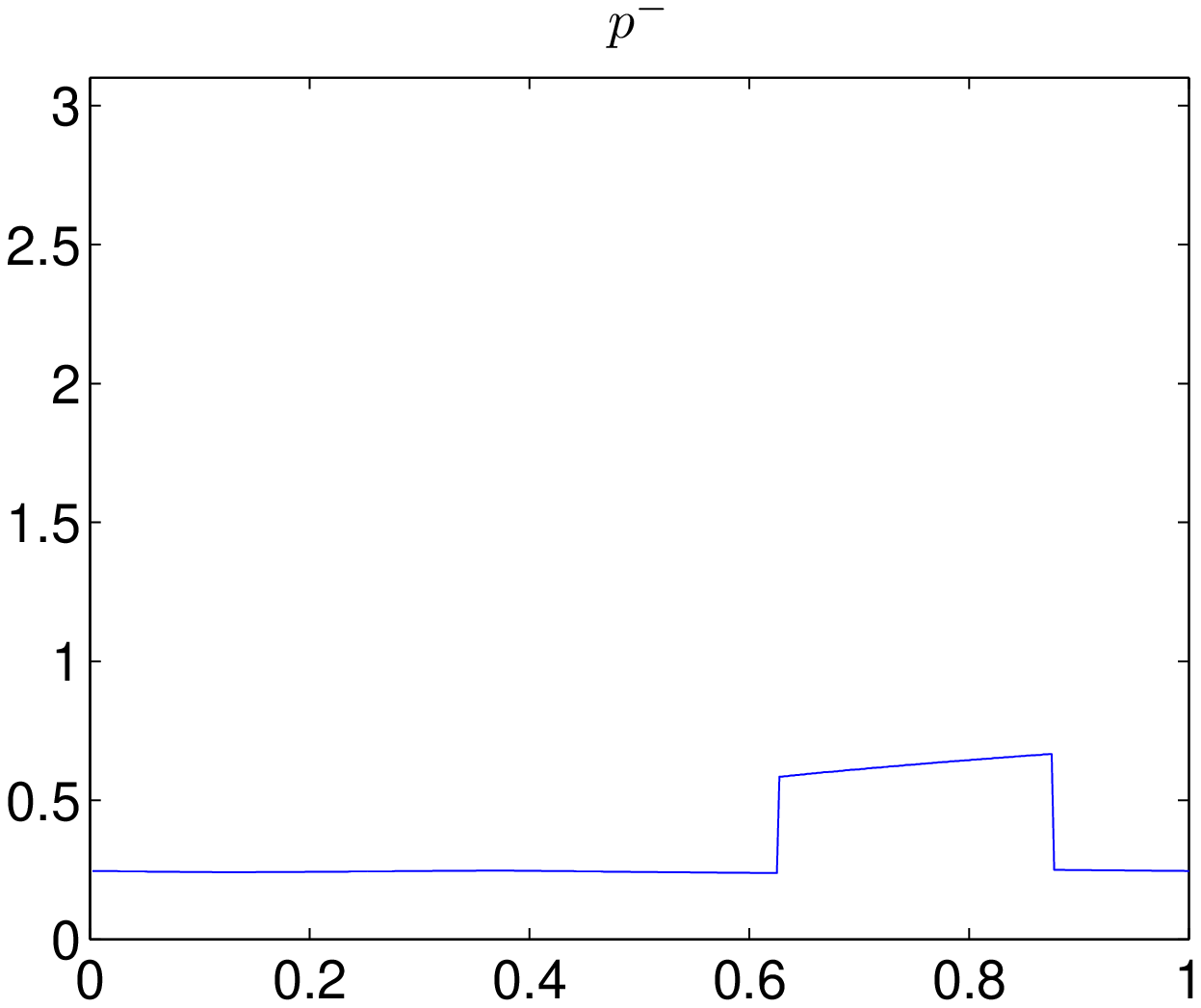}} &
\resizebox*{0.33\linewidth}{!}{\includegraphics{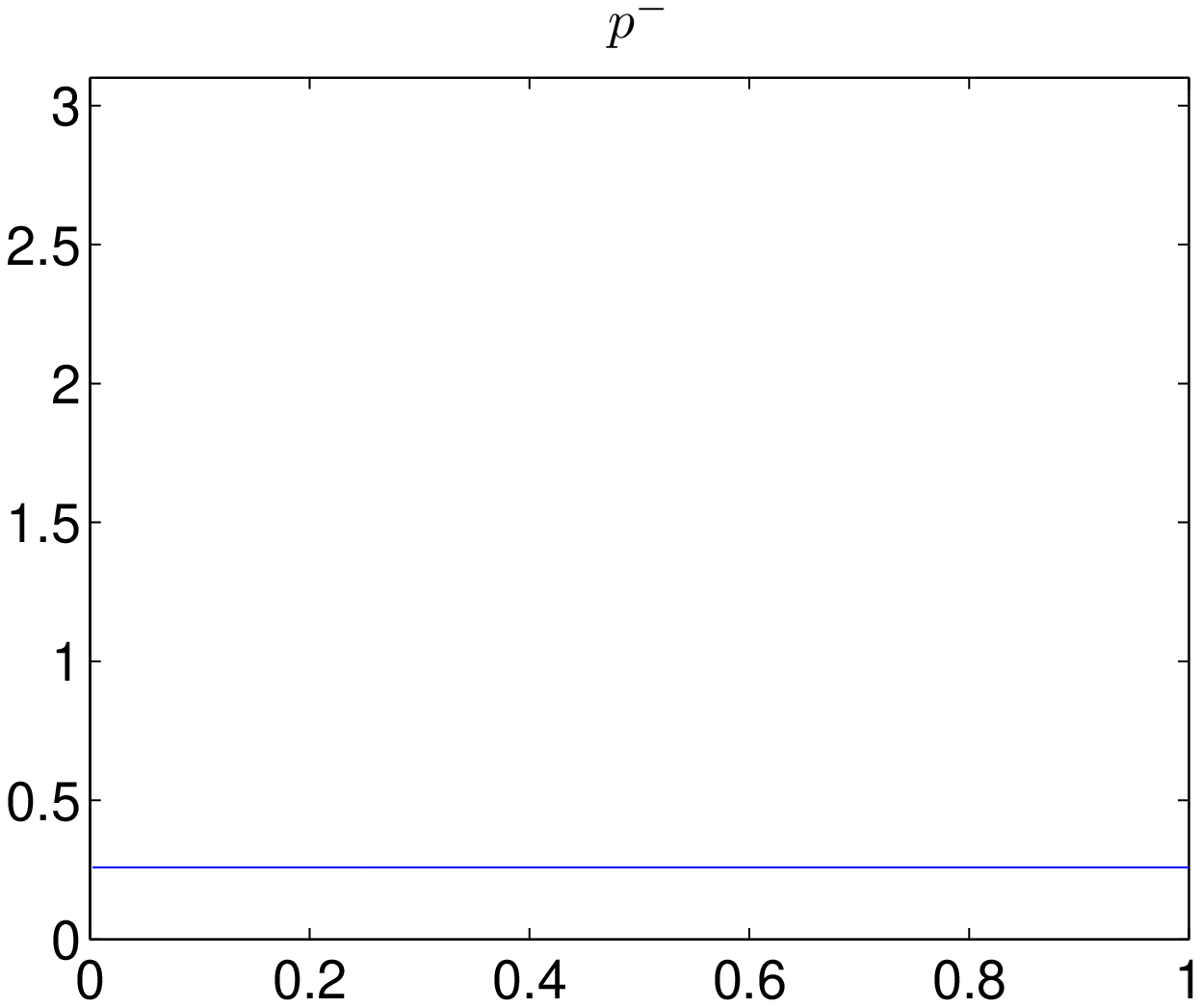}} \\
\resizebox*{0.33\linewidth}{!}{\includegraphics{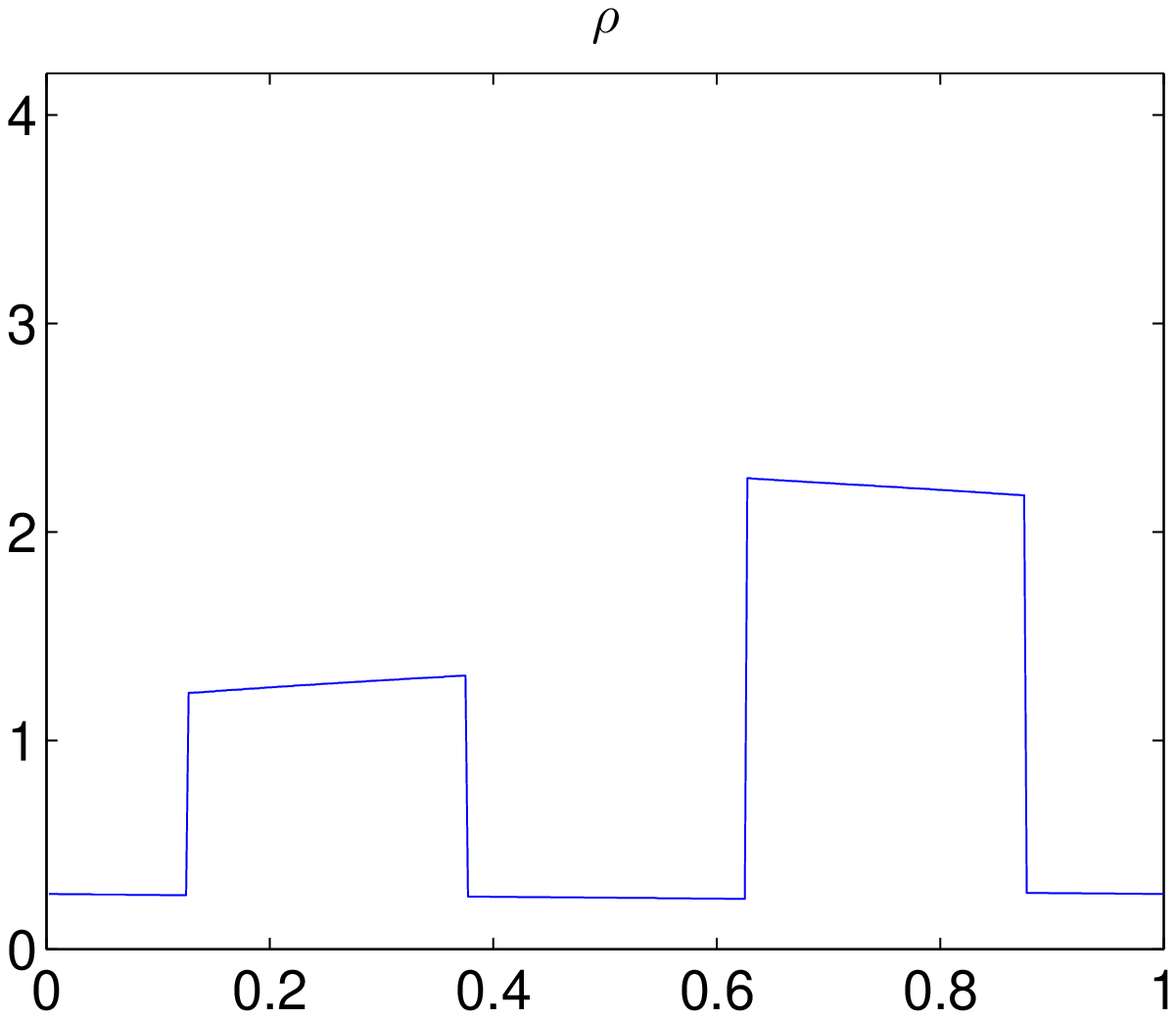}} &
\resizebox*{0.33\linewidth}{!}{\includegraphics{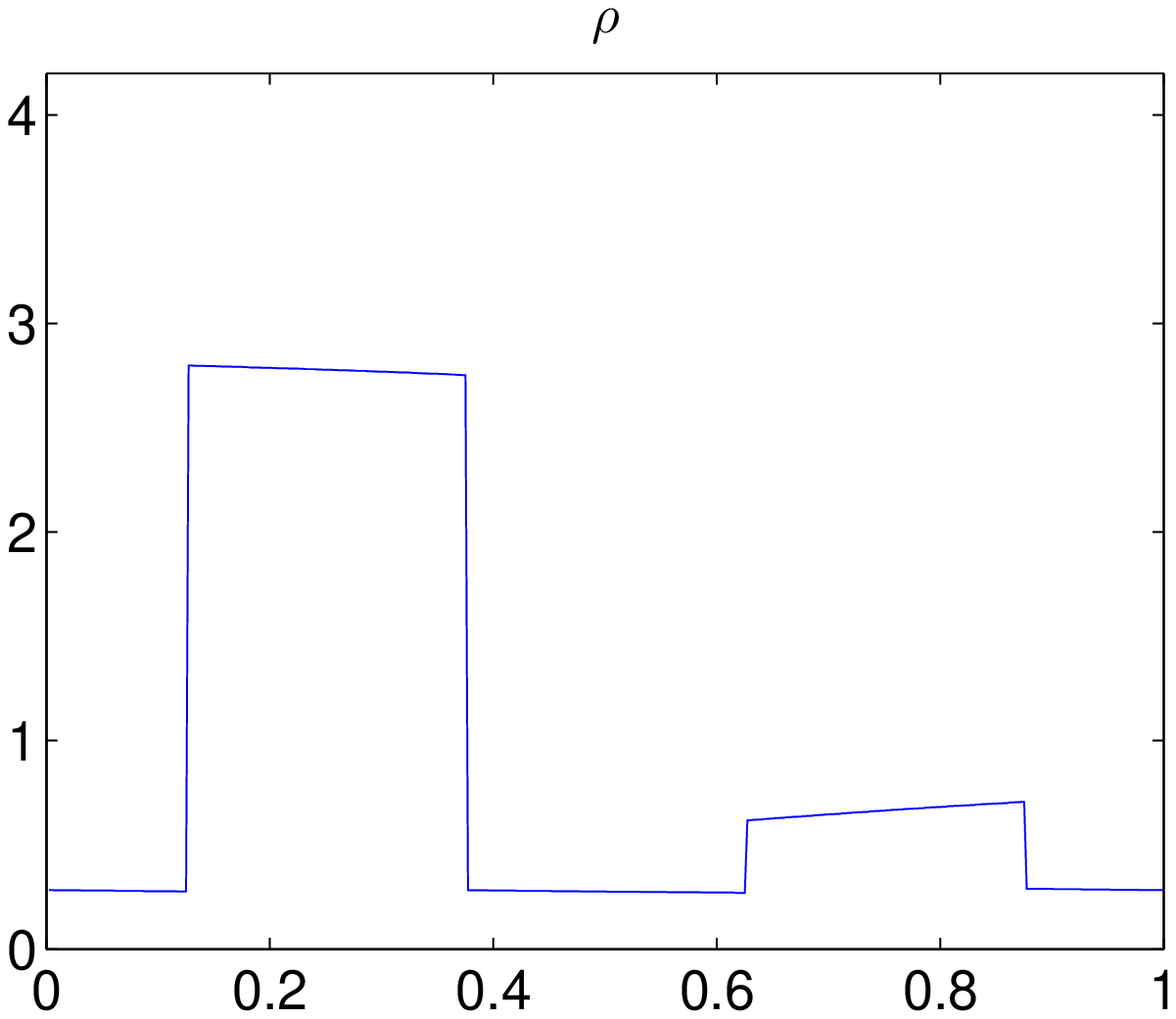}} &
\resizebox*{0.33\linewidth}{!}{\includegraphics{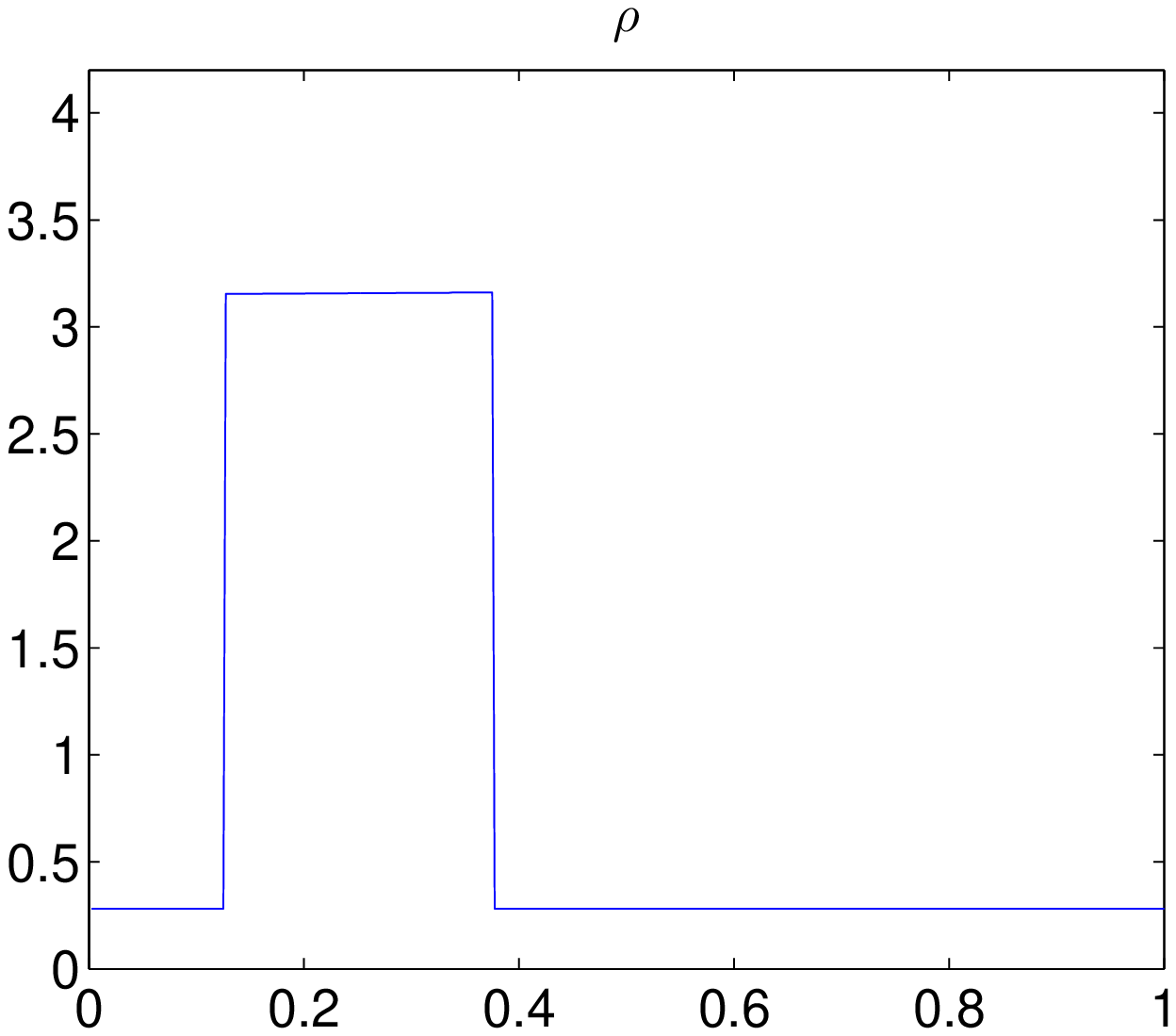}} \\
\resizebox*{0.33\linewidth}{!}{\includegraphics{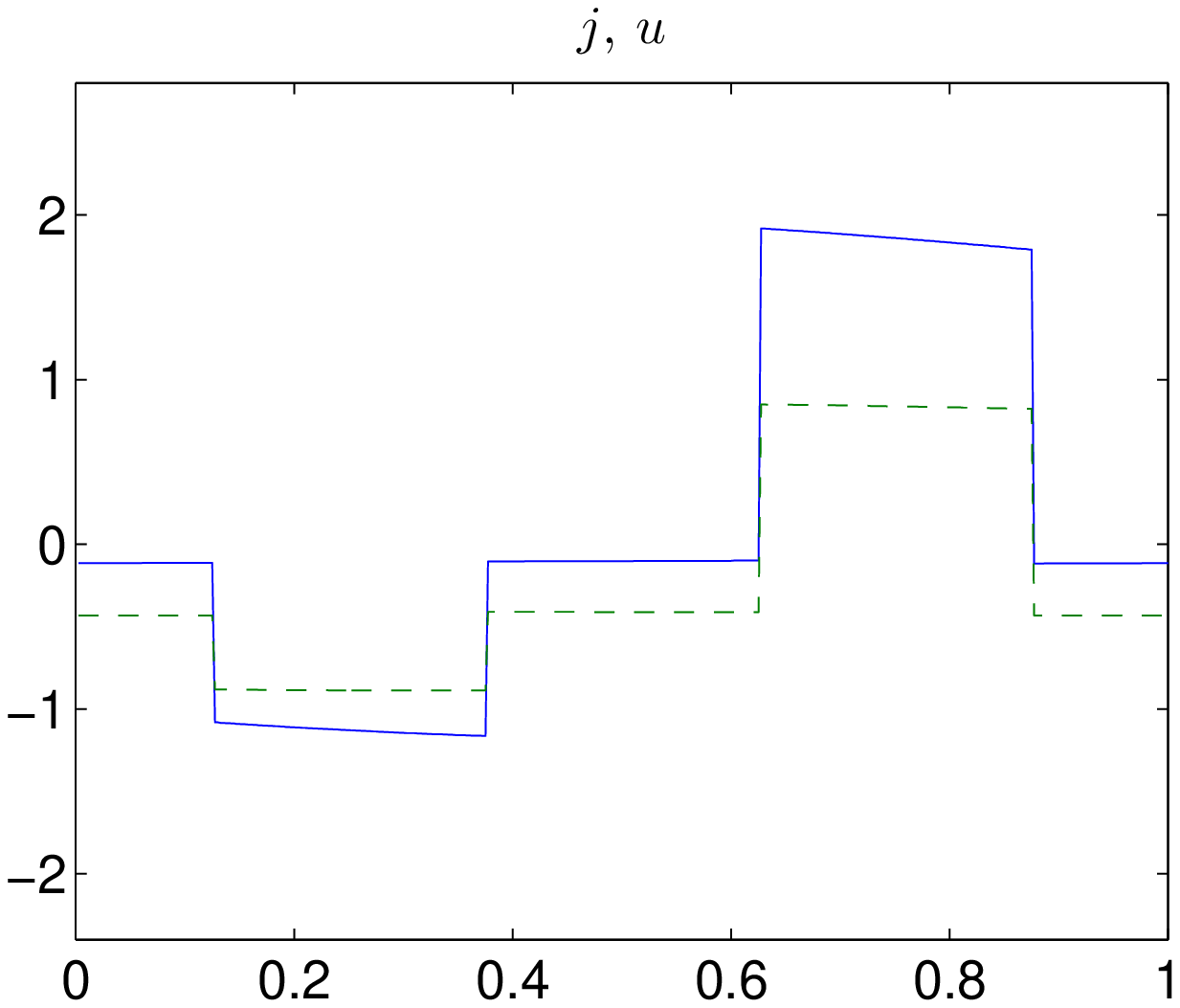}} &
\resizebox*{0.33\linewidth}{!}{\includegraphics{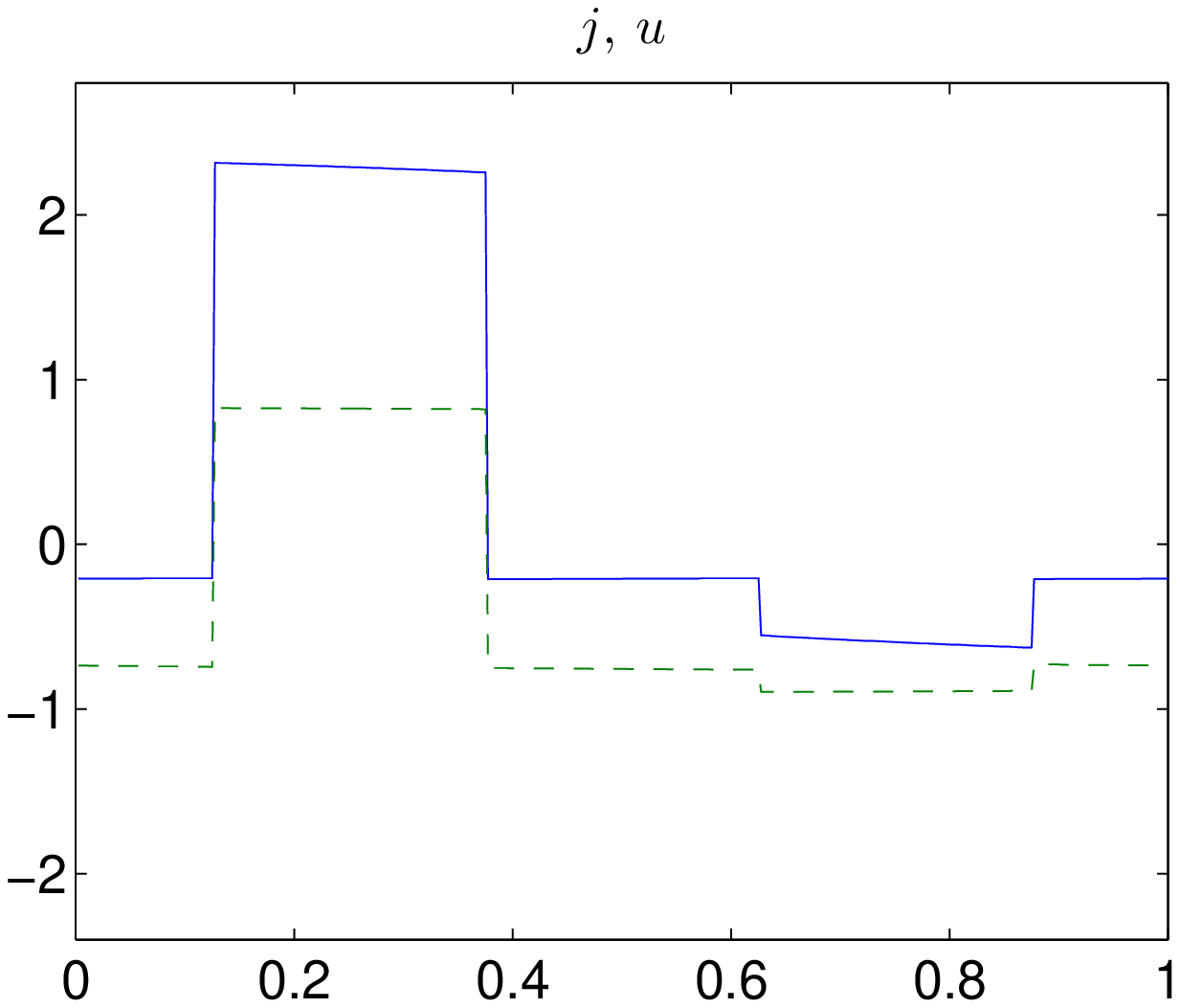}} &
\resizebox*{0.33\linewidth}{!}{\includegraphics{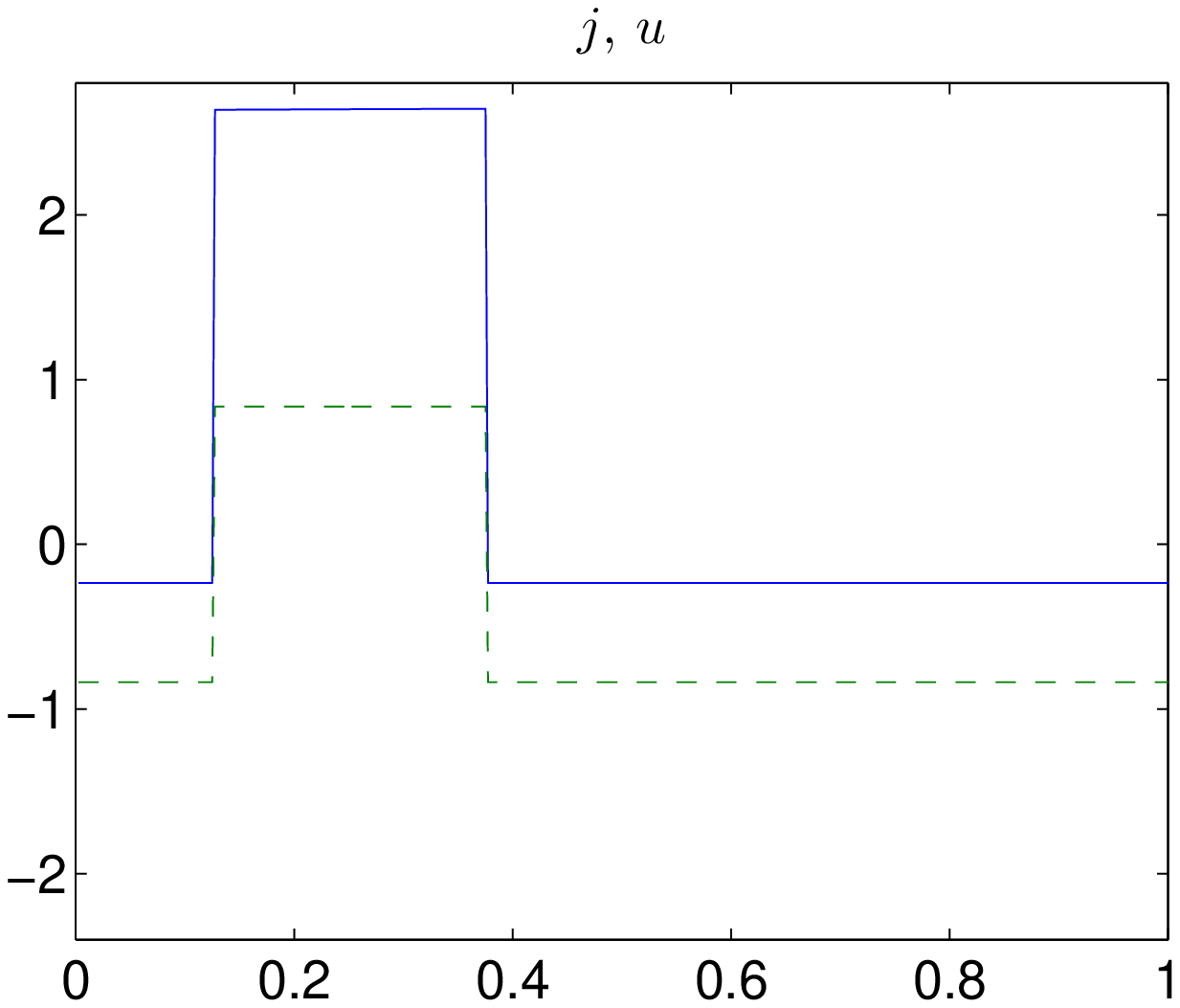}} \\
$t=2.5$ & $t=5.0$ & $t=15.0$
\end{tabular}\par}
\caption{Numerical results with the singular weight $\w=\delta_0$, $b=1$ and $\gamma_0=0.3$:
$p^-$ converges to a constant, while $p^+$ becomes a piecewise constant travelling wave profile
and $u$ (dashed line) jumps between the values $\pm u_s$ with $u_s = \sqrt{1-\gamma_0 b^{-1}} \cong 0.8367$.}
\end{figure}

In this subsection we briefly discuss the long time behaviour in the singular case $\w=\delta_0$.
Again, we start with a numerical example where we solve the kinetic system with the parameters $b=1$ and $\gamma_0=0.3$.
The initial condition is chosen as before, see Figure~\ref{fig:num1} (left panels).
In Figure~\ref{fig:num2} we present the time evolution of $p^+$\!\!, $p^-$\!\!, $\rho$, $j$ and $u$.
Based on the numerical observations, we conjecture that, for the small noise case $\gamma_0 < b$,
the long time dynamics are given by the travelling wave
profiles $(p_s^+, p_s^-)$, satisfying
\[
   \partial_t p_s^+ + \partial_x p_s^+ &=& 0 \,,\\
   \partial_t p_s^- - \partial_x p_s^- &=& 0 \,,
\]
and $\frac{p_s^+-p_s^-}{p_s^++p_s^-}\in \{u_s,-u_s\}$
with $u_s = \sqrt{1-\gamma_0 b^{-1}}$.
Then, it is a matter of a simple consideration to deduce
that one of the functions, say $p_s^-$, has to be a global constant,
while the other one, $p_s^+$, is a piecewise constant assuming only two values
$\{p^+_{s,1}, p^+_{s,2}\}$, satisfying the relations
\[
    p^+_{s,1} p^+_{s,2} = (p_s^-)^2 \qquad\mbox{and}\qquad
    \frac{p^+_{s,1}}{p^+_{s,2}} = \left(\frac{1-u_s}{1+u_s}\right)^2 \,.
\]
These relations chracterize the dynamic equilibria
between the densities of individuals marching to the left and to the right,
in dependence on the parameter values $\gamma_0$ and $b$.

\vspace{0.3cm}

For the large noise case ($\gamma_0 > b$), we prove that
the asymptotic state is again $j\equiv 0$ and $\rho\equiv 1$.
In fact, $j$ converges to $0$ in the $L^2$-sense exponentially fast as $t\to \infty$.
This follows from an application of the Gronwall lemma to the estimate
\[
    \frac12 \tot{}{t} \int_\Omega \left( \varrho^2 + j^2 \right) \d x
       &=& -2 \int (\gamma_0 + b(-1+u^2)) j^2 \d x \\
       &\leq& -2 (\gamma_0 - b) \int j^2 \d x \,.
\]

\section{The effect of shrinking interaction radius}\label{sec:Shrinking}
In biological applications, the weight function $\w$ is often chosen as the characteristic
function of the interval $[0,\sigma]$, where $\sigma>0$ is the interaction radius.
It is important to study the dependence of collective behaviour on the size of $\sigma$.
In this section, we will consider a theoretical limit $\sigma\to 0$.
Clearly, letting $\sigma\to 0$ with fixed $N$
leads to a trivial model without any interactions (it corresponds to the choice $b=0$).
On the other hand, letting first $N\to\infty$, followed by $\sigma\to 0$,
we obtain the model with $\w=\delta_0$, as mentioned in Section~\ref{sec:Existence}.
Consequently, the limits $\sigma\to 0$ and $N\to\infty$ do not commute.
From the modelling point of view, it is natural to study the limit $N\to\infty$
with the interaction radius shrinking as $N^{-\alpha}$ for some $\alpha>0$.
Indeed, as the space is getting more crowded, the visibility is reduced
and every individual can only take into account its closest neighbours.
Actually, taking $\w = \chi_{[0,N^{-\alpha}]}$ and letting $N\to\infty$,
we will show that a significant limit is obtained with $\alpha=1$,
leading to a new kinetic model.
The results are summarized in the following Theorem.
\begin{theorem}
Let $\w=\chi_{[0,N^{-\alpha}]}$. Then, the formal limit as $N\to\infty$
of the BBGKY hierarchy~$\fref{BBGKY}$ is
\begin{enumerate}[leftmargin = 1 cm, label = {(}\arabic*{)} ]
\item For $\alpha>1$:
\(
    \partial^2_t{\varrho} + 2(\gamma_0 + b)\,\partial_t{\varrho} = \partial^2_x{\varrho} \,.  \label{shrink1}
\)
\item For $0 < \alpha < 1$:
\(
    \partial_t{\varrho} + \partial_x{j} &=& 0\,,\label{shrink2}\\
    \partial_t{j} + \partial_x{\varrho} &=& -2 \left(\gamma_0 + b\left(1+\frac{j^2}{\varrho^2}\right)\right) j + 4bj \,.\label{shrink3}
\)
\item For $\alpha=1$:
\(
    \partial_t{\varrho} + \partial_x{j} &=& 0\,,\label{shrink4}\\
    \partial_t{\varrho} - \partial_x{j} &=& -2 \big(\gamma_0 + b(1+\eta)\big) j + 4b j \big(1-\exp(-\varrho)\big)  \,,\label{shrink5}
\)
with
\(
   \eta = \frac{j^2}{\varrho^2} \left( 1 - \exp(-\varrho) \right)
   + \left(1- \frac{j^2}{\varrho^2} \right) \exp(-\varrho) \left[ \mbox{Ei}(\varrho) - \gamma - \ln(\varrho) \right]  \,,  \label{eta}
\)
where $\mbox{Ei(z)}=\int_{-\infty}^z s^{-1} \exp(-s) \d s$ is the so-called exponential integral function
and $\gamma=0.577\dots$ is Euler's constant.
\end{enumerate}
\end{theorem}

\noindent
We see that for $0 < \alpha < 1$, we obtained the kinetic model with $u=j/\rho$
(i.e., the same as if one would first pass to $N\to\infty$ and then $\sigma\to 0$).
If $\alpha > 1$, we get the model with no interactions
(i.e., as if one would first pass to $\sigma\to 0$ and then $N\to\infty$)
which is described by the telegraph equation for $\varrho$ \cite{Kac:1974:SMR}.
In the significant limit with $\alpha=1$ we obtained a new model,
which is in fact the hydrodynamic model~\fref{hydrodynamic1}--\fref{hydrodynamic2}
with the function $\frac{j}{\varrho}(1-\exp(-\varrho))$ in place of $u$
and with $\eta$, given by \fref{eta}, in place of $u^2$.

\begproof
All we need to do is to recalculate the limits $N\to\infty$ in the expressions~\fref{q} and~\fref{r}
with $\w = \chi_{[0,N^{-\alpha}]}$.
Let us fix $z\in\Omega$ and assume that it is a Lebesgue point of $\varrho$ and $j$, i.e.,
\(   \label{Lebesgue}
    \lim_{\sigma\to 0} \frac{1}{\sigma} \int_0^\sigma \varrho(z-y) \d y = \varrho(z) \,,
\)
and similarly for $j$. Moreover, assume that $\varrho(z) > 0$.
In the same way as in Lemma~\ref{lemma:lim-q}, we calculate
\[
     Q_N &:=& \sum_{\boldv\in\setv{N-1}} \int_{\Omega^{N-1}}
          \frac{\sum_{m=1}^{N-1} \chi_{[0,N^{-\alpha}]}(|z-x_m|) v_m}{\sum_{m=1}^{N-1}
            \chi_{[0,N^{-\alpha}]}(|z-x_m|)} \prod_{i=1}^{N-1} p(x_i,v_i) \d\boldx \\
    &=& (N-1) \int_{\Omega^{N-1}}
           \frac{\chi_{[0,N^{-\alpha}]}(|z-x_1|)}{\sum_{m=1}^{N-1} \chi_{[0,N^{-\alpha}]}(|z-x_m|)}
                 \d P^{N-1}_\varrho(\boldx) \\
     &=& (N-1) K_N \int_{\Omega^{N-2}}
           \frac{1}{1+\sum_{m=1}^{N-2} \chi_{[0,N^{-\alpha}]}(|z-x_m|)} \d P^{N-2}_\varrho(\boldx) \,,
\]
where we denoted
\[
    I_N := \int_\Omega \chi_{[0,N^{-\alpha}]}(y) \varrho(z-y) \d y \,,\qquad
    K_N := \int_\Omega \chi_{[0,N^{-\alpha}]}(y) j(z-y) \d y \,.
\]
Obviously, the sum $\sum_{m=1}^{N-1} \chi_{[0,N^{-\alpha}]}(|z-x_m|)$ only takes the values $k=0, \dots, N-1$, and
\[
   P_\varrho^{N-2} \left( \left\{\boldx\in\Omega^{N-2};\, \sum_{m=1}^{N-2} \chi_{[0,N^{-\alpha}]}(|z-x_m|)
     = k\right\}\right) = \binom{N-2}{k} I_N^k (1-I_N)^{N-2-k} \,.
\]
Therefore,
\[
   \int_{\Omega^{N-2}} \frac{1}{1+\sum_{m=1}^{N-2} \chi_{[0,N^{-\alpha}]}(|z-x_m|)} \d P^{N-2}_\varrho(\boldx)
   &=& \sum_{k=0}^{N-2} \frac{1}{k+1} \binom{N-2}{k} I_N^k (1-I_N)^{N-2-k} \\
   &=& \frac{1}{(N-1) I_N} \left( 1 - (1-I_N)^{N-1} \right) \,.
\]
Due to~\fref{Lebesgue}, we have $\lim_{N\to\infty} N^\alpha I_N = \varrho(z)$
and $\lim_{N\to\infty} N^\alpha K_N = j(z)$.
Therefore,
\[
   \lim_{N\to\infty} Q_N = \lim_{N\to\infty} \frac{K_N}{I_N} \left( 1 - (1-I_N)^{N-1} \right)
        = \frac{j(z)}{\varrho(z)} \lim_{N\to\infty}  \left( 1 - (1-I_N)^{N-1} \right)\,,
\]
and
\[
    \lim_{N\to\infty} \left( 1 - (1-I_N)^{N-1} \right) =
                 \left\{ \begin{array}{cl}
                1-\exp(-\varrho(z)) & \mbox{for } \alpha=1 \,,\\
                1 &\mbox{for } 0 < \alpha < 1 \,,\\
                0 &\mbox{for } \alpha > 1 \,.
            \end{array}\right.
\]

Using the same technique as above, one calculates that
\[
    R_N &:=& \sum_{\boldv\in\setv{N-1}} \int_{\Omega^{N-1}}
          \left( \frac{\sum_{m=1}^{N-1} \chi_{[0,N^{-\alpha}]}(|z-x_m|) v_m}{\sum_{m=1}^{N-1} \chi_{[0,N^{-\alpha}]}(|z-x_m|)} \right)^{\!\!2} \,
          \prod_{i=1}^{N-1} p(x_i,v_i) \d\boldx \\
   &=& \frac{K_N^2}{I_N^2} \left( 1 - (1-I_N)^{N-1} \right) + \left(1 - \frac{K_N^2}{I_N^2}\right) G_N(I_N) \,,
\]
with
\[
   G_N (I_N) = \sum_{k=1}^{N-1} \frac{1}{k} \binom{N-1}{k} I_N^k (1-I_N)^{N-1-k} \,.
\]
It is easily shown that $G_N (I_N)$ vanishes in the limit $N\to\infty$ whenever $\alpha\neq 1$,
while with $\alpha=1$ the law of rare events (\cite{Feller:1967:IPT}) gives
\[
  \lim_{N\to\infty} G_N (I_N) = \exp(-\varrho(z)) \sum_{k=1}^\infty \frac{\varrho(z)^k}{k!\,k}
     = \exp(-\varrho(z)) \left[ \mbox{Ei}(\varrho(z)) - \gamma - \ln(\varrho(z)) \right] \,,
\]
We assume that $\varrho$ and $j$ be integrable functions,
so that almost every $z\in\Omega$ is a Lebesgue point for them.
\endproof

\begin{remark}\label{rem:shrink}
Our analysis can be seen as a first step towards a so-called topological model,
where the agent interactions are based on some connectivity graphs.
For instance, one can consider the situation where every agent
interacts only with its nearest neighbour. In this case, we have
a different time dependent weight function $\w_i$ for every agent, namely
$w_i = \chi_{[0,\sigma_i]}$ with $\sigma_i = \min_{m\neq i} |x_i-x_m|$.
Then, the passage to the corresponding continuum model
is a completely open problem; however, let us observe that
\[
    \left( \frac{1}{N-1} \sum_{m\neq i} \frac{1}{|x_i-x_m|^p} \right)^{-1/p} \mathop{-\!\!\!\longrightarrow}_{\!p\to\infty}\;
      \min_m |x_i-x_m| \,.
\]
Therefore, it would be interesting to consider the discrete system
with $\w_i = \chi_{[0,\sigma_i]}$,
\(   \label{sigma_i}
    \sigma_i = \left( \frac{1}{N-1} \sum_{m\neq i} \frac{1}{|x_i-x_m|^p} \right)^{-1/p} \,,
\)
and study the limit $N\to\infty$ and $p\to\infty$ (possibly with $p=N$).
Moreover, let us observe that ``on average'', $|x_i-x_m| \approx N^{-1}$. Therefore,
\[
    \left( \frac{1}{N-1} \sum_{m\neq i} \frac{1}{|x_i-x_m|^p} \right)^{-1/p} \approx N^{-1} \,,
\]
so we are in the situation of the significant limit $\alpha=1$ described above,
and we might believe that the new kinetic model obtained in this limit
can have some connection with the limit $N\to\infty$ and $p\to\infty$ of~$\fref{sigma_i}$.
\end{remark}

\section{Discussion}\label{Conclusions}
We introduced an individual based model with velocity jumps
aimed at explaining the experimentally observed collective motion of locusts
marching in a ring shaped arena~\cite{Buhl:2006:DOM}.
The frequency of individual velocity jumps increases with
a local or global loss of group alignment.
We showed that our model has the same predictive power
as the model of Czir\'ok and Vicsek, in particular,
it exhibits the rapid transition to highly aligned collective motion
as the size of the group grows and the switching of the group direction,
with frequency rapidly decreasing with increasing group size.
Moreover, in the limit $N\to\infty$ we obtained a system of two kinetic 
equations. We proved existence of its solutions and a partial result 
about the long time behaviour.
Finally, we studied the effect of shrinking the interaction radius $\sigma$
in the discrete model as the number of individuals, $N$, tends to infinity.
We showed that in the significant limit where $\sigma$ shrinks as $1/N$,
one obtains a new kinetic model.

Kinetic approach has previously been used in the literature to understand collective dynamics of individual based models.
Carrillo et al~\cite{Carrilo:2009:DMS} found the double milling phenomena
in the kinetic formulation of the model of self propeled particles with three zones of interactions.
The kinetic description of the Cucker-Smale model was introduced in~\cite{Ha:2008:PKH},
which can also be be derived from the Boltzmann-type equation, see~\cite{Carrilo:2009:AFD},
or Povzner-type equation,~\cite{Fornasier:2010:FDD}. For the survey of the most recent results
see the review~\cite{Carrillo:2010:PKH}.

Several interesting questions remain open, offering space for future investigations.
For example, the kinetic system~\fref{kinetic1}--\fref{kinetic2} deserves a better analysis,
in particular, uniqueness of solutions and more complete investigation of the long time behaviour.
It would also be interesting to know if and how the solutions corresponding to~$\w=\delta_0$
can be derived as a limit of solutions corresponding to $\w=\chi_{[0,\sigma]}$ as $\sigma\to 0$.
Another interesting direction of future research was formulated in Remark~\ref{rem:shrink}.
\vskip 3mm

\noindent{\bf Acknowledgment:} This publication was based 
on work supported by Award No. KUK-C1-013-04,
made by King Abdullah University of Science and Technology 
(KAUST).
JH acknowledges the financial support provided by 
the FWF project Y 432-N15 (START-Preis ``Sparse Approximation and Optimization in 
High Dimensions''). The research leading to these results 
has received funding from the European Research Council under
the {\it European Community's} Seventh Framework
Programme ({\it FP7/2007-2013})/ ERC {\it grant agreement} n$^o$ 
239870. RE would also like to thank Somerville College,
University of Oxford, for a Fulford Junior Research Fellowship.
Both authors would like to thank to the Isaac Newton Institute for
Mathematical Sciences in Cambridge (UK), where they worked together
during the program ``Partial Differential Equations in Kinetic Theories''.
The authors also acknowledge several interesting discussions
and valuable hints provided by Jan Vyb\'{\i}ral of the
Johann Radon Institute for Computational 
and Applied Mathematics (RICAM), Austrian Academy of Sciences,
and Christian Schmeiser of the Faculty of Mathematics, University of Vienna.


{\small

\bibliographystyle{siam}
\bibliography{bibrad}

}

\end{document}